\documentclass[preprint,review,12pt]{elsarticle}
\usepackage[margin=2.5cm]{geometry}

\usepackage{amssymb}
\usepackage{color}
\usepackage{amsthm}
\usepackage{amsmath}
\usepackage{amscd,amsfonts,bbold}
\usepackage{caption}
\usepackage{subcaption}
\usepackage{color,soul}
\usepackage{lineno}
\usepackage{tikz-cd}
\usepackage{array}
\newcolumntype{P}[1]{>{\centering\arraybackslash}p{#1}}
\newcolumntype{M}[1]{>{\centering\arraybackslash}m{#1}}
\usepackage{hyperref}

\journal{International Journal for Computational Methods in Engineering Science and Mechanics, September 9, 2017}

\begin{document}

\begin{frontmatter}



\title{Numerical Convergence of Discrete Exterior Calculus on Arbitrary
	Surface Meshes}


\author[address1]{Mamdouh S. Mohamed\corref{cor1}\fnref{fn01}}
\ead{mamdouh.mohamed@kaust.edu.sa}
\author[address2]{Anil N. Hirani}
\ead{hirani@illinois.edu}
\author[address1]{Ravi Samtaney}
\ead{ravi.samtaney@kaust.edu.sa}

\address[address1]{Mechanical Engineering, Physical Science and Engineering Division, King Abdullah University of Science and Technology, Thuwal, KSA}
\address[address2]{Department of Mathematics, University of Illinois at Urbana-Champaign, IL, USA}

\cortext[cor1]{Corresponding author}
\fntext[fn01]{On leave, Department of Mechanical Design and Production, Faculty of Engineering, Cairo
University, Giza, Egypt.}

\begin{abstract}

Discrete exterior calculus (DEC) is a structure-preserving numerical framework for partial differential equations solution, particularly suitable for simplicial meshes. A longstanding and widespread assumption has been that DEC requires special (Delaunay) triangulations, which complicated the mesh generation process especially on curved surfaces. This paper presents numerical evidences demonstrating that this restriction is unnecessary. Convergence experiments are carried out for various physical problems using both Delaunay and non-Delaunay triangulations. Signed diagonal definition for the key DEC operator (Hodge star) is adopted. The errors converge as expected for all considered meshes and experiments. This relieves the DEC paradigm from unnecessary triangulation limitation.  
    
\end{abstract}


\end{frontmatter}


\section{Introduction}
\label{sec:introduction}

The discrete exterior calculus (DEC) framework \cite{hirani2003discrete,desbrun2005discrete} provides discrete definitions for many exterior calculus operators, aiming mainly to numerically solve systems of equations arising in field theories, such as fluid dynamics and electromagnetism, and is also used in computer graphics applications. The main characteristics of DEC operators are their mimetic behavior, retaining at the discrete level many of the identities/rules of their smooth counterparts, and their coordinate-independent nature. The mimetic character of DEC operators provides conservation properties in the resulting numerical methods. For example, in DEC discretization of incompressible Navier-Stokes equations \cite{mohamed2016discrete}, both mass and vorticity conservation are a direct result of this mimetic character of DEC operators. Some properties, such as kinetic energy conservation for Euler equation, may require in addition the construction of a conservative time discretization as well. DEC was applied during the past decade to numerically solve various physical problems including Darcy \cite{HiNaCh2015,gillette2011dual} and incompressible Navier Stokes flows \cite{elcott2007stable,mullen2009energy,mohamed2016discrete}. 

DEC is one of many mimetic discretizations that have appeared in literature. A related discretization of exterior calculus is finite element exterior calculus (FEEC) \cite{arnold2010finite} which differs from DEC in the way that metric information is treated. While FEEC permits arbitrary order accuracy, DEC is limited to lowest order. Some of the other mimetic discretizations include the covolume method \cite{HallCavendish:1991dvmsvf,cavendish1992solution,nicolaides:1992direct,NiTr2006} and those based on the support operator method \cite{HySh1997a,HySh1997b}. The covolume method is defined for planar and three-dimensional domains on which it is identical to DEC. The support operator based methods were originally defined also for such domains and can be considered as a discretization of classical vector calculus. DEC may be characterized as a generalization of older mimetic methods such as the covolume method. The separation of metric and non-metric operators in exterior calculus and the implementation of this separation in DEC and FEEC make DEC/FEEC suitable for piecewise linear manifolds of any dimension, with a useful case being surfaces in $\mathbb{R}^3$. 

A characteristic requirement of DEC is the usage of a dual mesh in addition to the primal simplicial mesh. The dual mesh enables the DEC operators to discretely mimic many key calculus theorems; e.g. Stokes theorem and divergence theorem. The common choice for a dual mesh in DEC implementations is the circumcentric dual. This is mainly due to the mutual orthogonality between primal and dual mesh objects, resulting in diagonal definitions for all Hodge star operators and their inverses. It has been a commonly held belief that one of the limitations of the circumcentric dual is that it is well-defined only on Delaunay triangulations. For non-Delaunay triangulations, the direction of some circumcentric dual edges is flipped with respect to the considered conventional orientation. Moreover, the dual cells associated with some primal nodes end up overlapping. Such complications led to the notion that using a circumcentric dual on a non-Delaunay mesh might result in incorrect DEC solutions.
 
This idea that circumcentric duals may not be used on non-Delaunay triangulations dates back to the covolume method literature. The covolume method employs operators definitions similar to the DEC operators albeit with a different terminology. The method was used mainly to solve the incompressible Navier-Stokes equations both on 2D flat surfaces and in 3D domains \cite{nicolaides:1989flow,choudhury1990discretization,HallCavendish:1991dvmsvf,cavendish1992solution,hu1992covolume,hall1992network,nicolaides:1992direct,nicolaides1993covolume,cavendish1994complementary}. It was clearly stated in most covolume literature that such a scheme is developed for Delaunay/Voronoi pairs, representing the primal/dual meshes. There exist, however, a remark by Nicolaides \cite{nicolaides1993covolume} that the covolume method might work on non-Delaunay meshes with additional coding effort to deal with the non-simple (overlapping) dual cells. Nevertheless, no further details or tests were provided to further investigate or verify such a remark. In the work by Perot and coauthors \cite{perot2000conservation,zhang2002accuracy} on the other hand, it was pointed out that the covolume algorithm can still hold in principle for arbitrary meshes but with the possibility of being intolerably inaccurate for such arbitrary meshes. The application of the developed covolume algorithms was consequently limited only to Delaunay meshes.

Avoiding the circumcentric dual on non-Delaunay triangulations became the common practice in discrete differential geometry and exterior calculus implementations. For example, when defining discrete differential geometry operators; e.g. normal vector and curvatures, Meyer et al. \cite{meyer2002discrete} reverted locally to the barycentric dual mesh whenever a non-Delaunay mesh pair exists. Moreover, DEC discretizations of the incompressible Navier-Stokes equations \cite{mullen2009energy,mohamed2016discrete} employed the circumcentric dual mesh and consequently considered only Delaunay triangulations for the application of the developed schemes. Furthermore, Mullen et al. \cite{mullen2011hot} pointed out that for a Delaunay/Voronoi triangulation, failure to keep the circumcenter inside its triangle/tetrahedron can lead to numerical degeneracy. Hirani et al. \cite{hirani2013delaunay} showed that for a certain choice of sign convention for circumcentric dual objects, the discrete Hodge star assembled based on such elementary dual pieces is positive definite if the primal mesh is Delaunay. Delaunay meshes are therefore sufficient for obtaining positive entries in the diagonal Hodge matrices. In that paper, no mathematical claim was made regarding the solution accuracy in case the Delaunay condition was violated. However, a numerical experiment in Hirani et al. \cite{hirani2013delaunay} showed the solution of a scalar Poisson equation on a non-Delaunay triangulation (using the circumcentric dual) to be incorrect. Later on, an error in the code was discovered \cite{HiKaVa2017}. After correcting the code, the solution on both Delaunay and non-Delaunay meshes appeared to be similar and correct.

The literature in general leaned towards the expectation that using circumcentric duals on non-Delaunay triangulations can lead to incorrect results. However this claim remained largely untested and was thought to be one of the main limitations of both the covolume method and the DEC scheme. For applications involving non-Delaunay triangulations, the barycentric dual was alternatively chosen. This required alternative definitions, like the Galerkin \cite{bossavit1998computational} and the barycentric \cite{trevisan2004geometric,auchmann2006geometrically} definitions, for some of the Hodge star operators. The reader may refer to \cite{mohamed2016baryHodge} for numerical DEC implementations using such operators. Both the Galerkin and the barycentric definitions have, however, a sparse non-diagonal representation, which prohibits a sparse representation for their inverse operators. This is in addition to the increased computational cost when using the Galerkin or the barycentric definitions Hodge star definitions \cite{mohamed2016baryHodge}, in comparison to the diagonal definitions based on the circumcentric duals. Such  alternative dual meshes, alternative Hodge star definitions and the accompanied complications may not in fact be necessary if one were able to demonstrate that DEC solutions using the circumcentric dual work correctly on non-Delaunay triangulations. Motivated by increasing the range of applicability of both the DEC and the covolume methods, an investigation  of utilizing circumcentric duals on non-Delaunay meshes is warranted. This is the focus of the present work. 

The outline of this paper is as follows. A brief overview of the convention used in constructing the circumcentric dual and defining the Hodge star operators on arbitrary meshes is first discussed in Section \ref{sec:DEC-arbitrary}. Section \ref{sec:results} presents several numerical experiments for the solution of the 0-form Poisson, the 1-form Poisson and the incompressible Navier-Stokes equations on non-Delaunay meshes. Triangulations having various proportions of non-Delaunay triangle pairs are considered during the current study. The results also explore how the stiffness matrices condition number changes when considering non-Delaunay meshes. The paper closes with concluding remarks emphasizing the key findings and the significance of the presented research.          

\section{The circumcentric dual on arbitrary triangulations}
\label{sec:DEC-arbitrary} 

The current research focuses only on simplicial meshes approximating flat/curved surfaces. The simulation domain $\Omega$ is therefore considered to have a dimension $N=2$. The constituting elements of the simplicial complex $K$ are primal $k$-simplices $\sigma^k \in K$, each is defined by the indices of the nodes forming it as $\sigma^k = [v_0, ..., v_k], k \le 2$. The dual to a primal $k$-simplex $\sigma^k$ is the $(N-k)$-cell denoted by $\star \sigma^k \in  \star K$, where $\star K$ is the dual mesh complex. The orientation convention for both the primal triangles and the dual cells is counterclockwise. For primal edges shared by Delaunay triangles, the orientation convention of their dual edges is $90^\circ$ counterclockwise with respect to primal edges orientation.      

\begin{figure}
	\centering
	\includegraphics[width=0.49\textwidth]{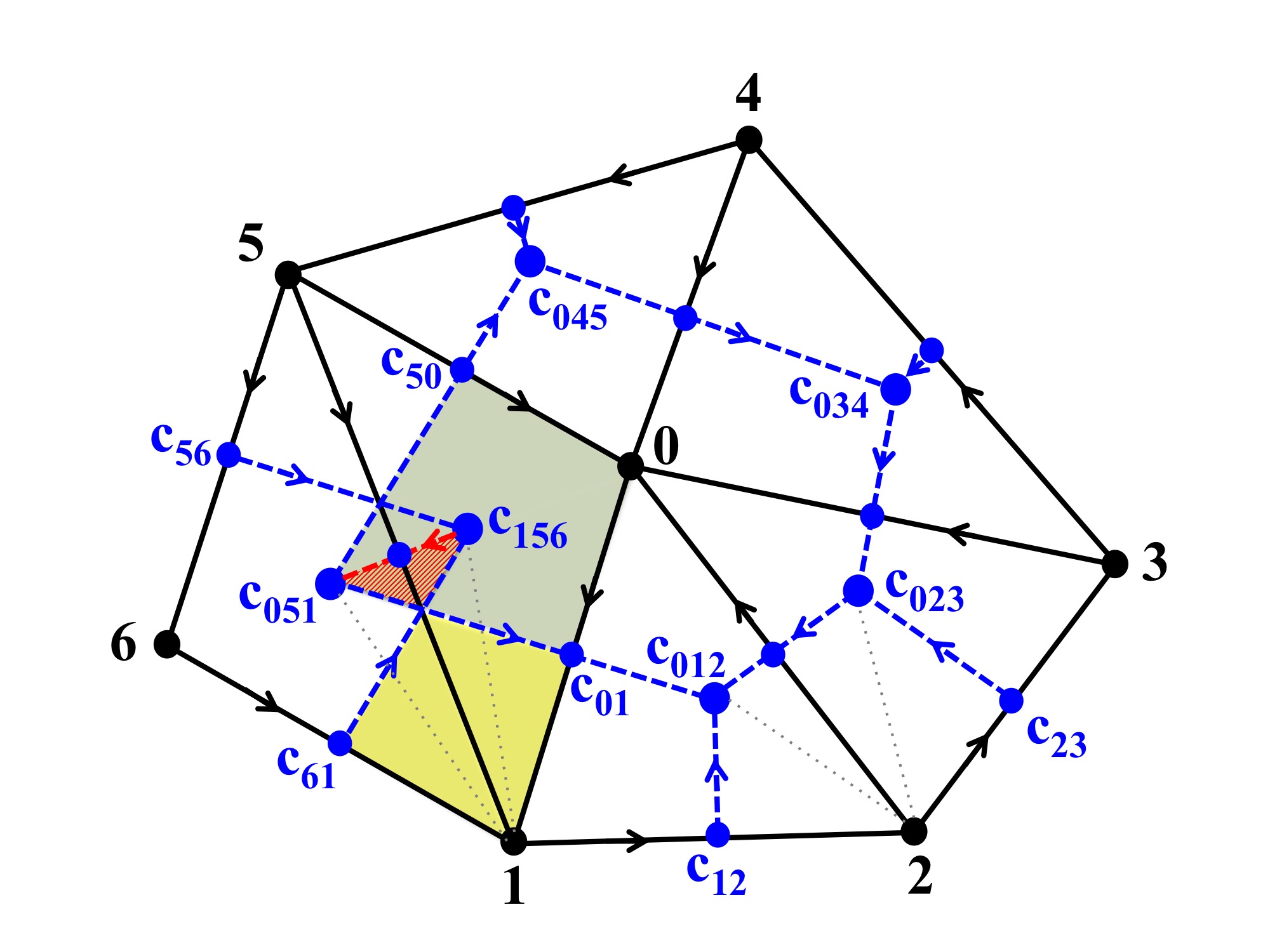}
	\caption{A sketch of a sample non-Delaunay mesh showing the primal edges in solid black color and the circumcentric dual edges in dashed blue color. The dual edge shared by the non-Delaunay triangles pair is shown in dashed red color. The dual area associated with node $[0]$ (only the contribution from triangle $[0,5,1]$) is shown in green color and has a positive area. The dual area associated with node $[1]$ (the contribution from both triangles $[0,5,1]$ and $[1,5,6]$) consists of the yellow sector (having a positive area) and the shaded red sector (having negative area).}
	\label{Fig:mesh}
\end{figure}

We consider arbitrary simplicial triangulations during the current study. Figure \ref{Fig:mesh} shows a sample mesh that includes a non-Delaunay triangles pair $[0,5,1]$ and $[1,5,6]$. The process of constructing the circumcentric dual commences by defining the nodes dual to the primal triangles as their circumcenters. This is followed by connecting the nodes dual to neighbor triangles through the shared primal edge midpoint to define the dual edges. Finally, the area dual to a primal node is defined as the area enclosed by the edges dual to the primal edges sharing this node. The circumcentric dual mesh tiles the whole domain without any cells overlapping only if the primal simplicial mesh is Delaunay. For non-Delaunay triangulations, some dual cells overlap with each other and some dual edges have their conventional direction flipped. For the purpose of defining the Hodge star operators, it is important to correctly consider the \emph{signed} volume (i.e., length or area) for these flipped dual edges and overlapping dual cells. Such a convention for the signed volumes was discussed in more detail by Hirani et al. \cite{hirani2013delaunay}. Relevant details are briefly presented here for completeness. 

The primal edge $[5,1]$ is shared by the non-Delaunay pair $[0,5,1]$ and $[1,5,6]$. According to edge $[5,1]$ orientation with respect to both triangles, the dual edge $\star[5,1]$ is defined to point from $c_{156}$ to $c_{051}$, where $c_{ijk} = \star[i,j,k]$ is the node dual to the triangle $[i,j,k]$. Such an orientation would make $\star[5,1]$ to be oriented $90^\circ$ counterclockwise with respect to $[5,1]$ orientation in case the neighbor triangles were Delaunay. However, since triangles $[0,5,1]$ and $[1,5,6]$ are non-Delaunay, the dual edge $\star[5,1]$ (shown in dashed red color in Fig. \ref{Fig:mesh}) has a flipped direction and is oriented $90^\circ$ clockwise with respect to edge $[5,1]$ orientation. According to the volume convention in \cite{hirani2013delaunay}, the length of such flipped dual edge is considered to be negative.

The negative sign attached to flipped dual edges is not only considered when calculating their own volumes, but further when calculating their contributions to the dual cells volumes. For the primal node $[2]$ in Fig. \ref{Fig:mesh}, for example, its dual cell area is calculated as the summation of the area sectors formed by node $[2]$ and the dual edges $[c_{23},c_{023}]$, $[c_{023},c_{012}]$ and $[c_{12},c_{012}]$. Since none of the dual edges is flipped, the area of all sectors is considered to be positive. Following the same methodology for node $[1]$, the dual cell $\star[1]$ (considering for the moment only the contribution from both triangles $[0,5,1]$ and $[1,5,6]$) consists of the area sectors formed by node $[1]$ and the dual edges $[c_{61},c_{156}]$, $[c_{156},c_{051}]$ and $[c_{051},c_{01}]$. Since the dual edge $[c_{156},c_{051}]$ is flipped, its contribution to $\star[1]$ is considered to be negative. This results in having the dual cell $\star[1]$ to consist of a positive volume sector (colored in yellow in Fig. \ref{Fig:mesh}) and a negative-volume sector (shaded in red in Fig. \ref{Fig:mesh}), where the net dual cell area is the signed summation of both. The negative-volume sector from the dual cell $\star[1]$ overlaps with the positive sectors from both the dual cells $\star[0]$ and $\star[6]$. The dual cell $\star[5]$ is similar to $\star[1]$ in both having a negative-volume sector and its overlapping with $\star[0]$ and $\star[6]$. It is worth noting that following such sign convention makes a negative-volume sector to overlap with two positive volume sectors. This results in having the total domain area calculated by either  summing the primal triangles areas or the dual cells areas to ultimately coincide. Accordingly, for non-Delaunay triangulations, the circumcentric dual cells tile the whole domain in a signed sense. Whenever an integration is carried out over the dual cells, considering the sign of these negative-volume sectors will make each point in the domain to be considered only once eventually. Such dual cell integration is carried out in the flux error calculation during the incompressible Navier-Stokes solution in Section \ref{subsec:Incomp-Navier-Stokes}.  

The negative signs attached to the volumes of flipped dual edges and some dual cell sectors, in the case of non-Delaunay triangulations, can be further justified in a physical sense. Consider an incompressible fluid flow given by a stream function 0-form $\psi$ defined on the primal nodes. The velocity 1-form $u$ defined on the dual edges is then calculated as $u=\ast_1 \textrm{d}_0 \psi$, representing the velocity vector integration along a dual edge. The primal 1-form $\textrm{d}_0 \psi$ represents the fluid mass flux across a primal edge. Multiplying this mass flux by the Hodge star operator $\ast_1$ (dividing by the primal edge length and multiplying by the dual edge length) provides a first order approximation of the velocity integration along the dual edge; i.e. the dual 1-form $u$. Whenever a primal edge is shared by non-Delaunay triangles pair, the 1-form $\textrm{d}_0 \psi$ continues to represent the mass flux, whereas the flipped dual edge points opposite to its conventional direction. Therefore, a negative sign is required to be added to the Hodge star operator $\ast_1$ in order to account for the dual edge flipped direction. This negative sign is added through the negative volume considered for such flipped dual edge.  

In regards to the dual cells, we consider calculating the vorticity 0-form on the primal nodes as $\omega = \ast_0^{-1}  [-\textrm{d}_0^T] u$. The dual 2-form $[-\textrm{d}_0^T] u$ calculates the circulation integration on the dual cells in a way that discretely mimics Stokes theorem through integrating the velocity field on the dual cells boundary. For the dual cells whose boundaries do not include flipped dual edges, the boundary integration is carried out in a counterclockwise sense; e.g. the dual cell $\star [0]$ in Fig. \ref{Fig:mesh}. Otherwise, the boundary integration is carried out following a clockwise convention on the negative-volume dual cell sectors. Therefore, for the dual cell $\star [1]$ in Fig. \ref{Fig:mesh}, the action of the $[-\textrm{d}_0^T]$ operator, mimicking a boundary integration process, is carried out in a counterclockwise sense over the yellow positive-volume sector, but in the opposite sense over the red negative-volume sector. Accordingly, when calculating a nodal averaged vorticity value through multiplying by the Hodge star operator $\ast_0^{-1}$, it is required to divide the circulation integration by the net dual cell area calculated as a summation of positive and negative area sectors in order to account for such opposite conventions during the boundary integration process. This can be better perceived by considering the limit where the area integrated quantity is almost constant. Since the integration is carried out over various sectors following opposite conventions, this is equivalent to an integration over the net signed area following a convention consistent with the net area sign. Consequently, averaging the integrated quantity should be through dividing by the signed net area.

The construction of such circumcentric dual mesh enables diagonal definitions for all Hodge operators, and also their inverses, as $\ast_k = \frac{|\star \sigma^k|}{|\sigma^k|}$ for $k=0,1,2$, where $|.|$ is the enclosed simplex/cell signed volume. Whenever non-Delaunay triangle pairs exist, the flipped dual edges volume is considered to be negative. Similarly, the area of some dual cells can, in principle, have a net negative value. During the experiments presented in the next section, some non-Delaunay meshes did have dual cells with a net negative area. The ability of such signed diagonal definitions for the Hodge star operators to provide correct numerical solutions is investigated in the next section.                     
 
\section{Results and discussion}
\label{sec:results}

Convergence tests are performed for four physical test cases aiming to benchmark the behavior of the circumcentric diagonal Hodge star operators defined on non-Delaunay triangulations. The four test cases are the primal 0-form Poisson equation, the dual 0-form Poisson equation, the primal 1-form Poisson equation and the incompressible Navier-Stokes equations. For all Poisson equations, we use a standard formulation rather than a mixed formulation \cite{arnold2010finite}. The mixed formulation involves a secondary unknown and eliminates the need for inverse Hodge star operators. However, this is not necessary during DEC solutions since a diagonal representation for the inverse Hodge star operators is available.

The convergence tests for each of the physical test problems described above are carried out using four groups of meshes. The first group includes five Delaunay meshes created independently, where the number of triangular elements in a finer mesh is almost four times the number of triangles in the coarser mesh. The remaining three groups are non-Delaunay triangulations having various proportions of non-Delaunay triangles pairs, and are generated by distorting the Delaunay triangulations. The distortion process is carried out as follows: for a given Delaunay triangulation, a random primal edge is selected and the apexes of the two neighbor triangles are moved towards the edge's midpoint to make these two triangles a non-Delaunay pair. The distortion process continues until a specified ratio of edges shared by non-Delaunay triangle pairs is achieved. 

\begin{figure}
		\centering
		\includegraphics[width=0.49\textwidth]{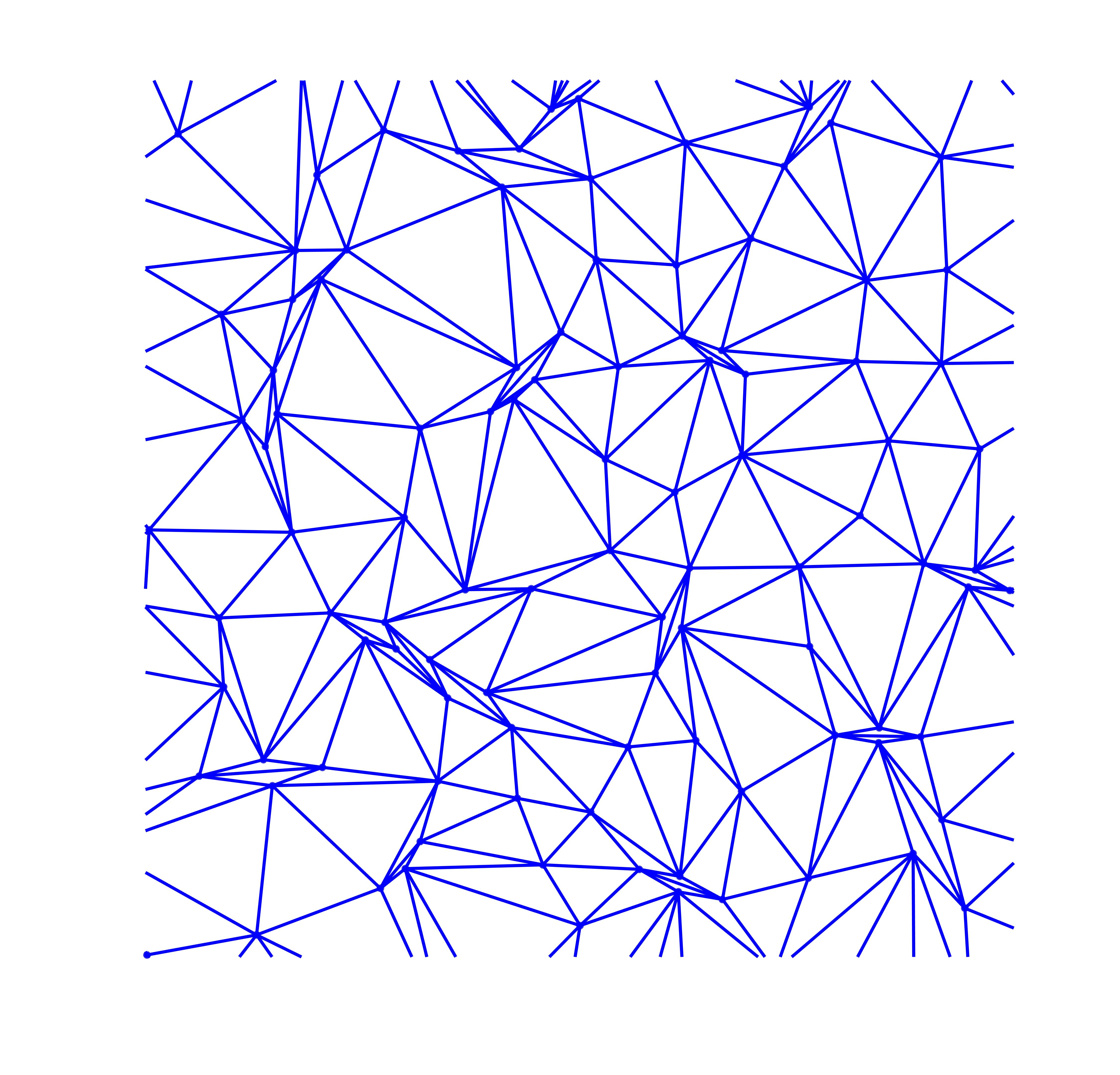}
		\caption{Part of a sample mesh corresponding to the second coarsest mesh from the non-Delaunay mesh group 3.}
		\label{Fig:meshSample}
	\end{figure}

For the three groups of non-Delaunay triangulations, the ratios of edges shared by non-Delaunay triangles are almost $1\%$, $5\%$ and $15\%$, as stated in Table \ref{table:mesh-types}. The table also provides the ratio of the non-Delaunay triangles: this is calculated by flagging both triangles in each non-Delaunay pair and then calculating the ratio of flagged triangles to the total number of triangles. The triangles maximum aspect ratio values indicate the non-uniformity level of the considered meshes. For some non-Delaunay meshes, many of the triangles with a boundary edge have their circumcenter residing outside the whole simulation domain. A sample non-Delaunay triangulation is shown in Fig. \ref{Fig:meshSample}.

	\begin{table}[h!]
		\centering
		\begin{tabular}{|M{3.2cm}|M{3cm}|M{3cm}|M{1.9cm}|}
			\hline
			& Edges having non-Delaunay triangles pair & Non-Delaunay triangles ratio & Max aspect ratio\\
			\hline
			Non-Delaunay mesh 1  & $\sim 1\%$& $\sim 3\%$ & $76$ \\
			\hline
			Non-Delaunay mesh 2  & $\sim 5\%$& $\sim 15\%$ & $97$ \\
			\hline
			Non-Delaunay mesh 3  & $\sim 15\%$& $\sim 43\%$ & $288$ \\
			\hline
		\end{tabular}
		\caption{The groups of non-Delaunay meshes used during the numerical solutions. The table shows the ratio of edges shared by non-Delaunay triangles pairs, the ratio of non-Delaunay triangles and the maximum aspect ratio for the triangles (the circumradius divided by the inradius). The numbers stated are for the finest mesh from each mesh group, however they are representative to all meshes in the group.}
		\label{table:mesh-types}
	\end{table}
	
The conducted numerical experiments demonstrate the $L^2$ error convergence, characterizing the difference between the numerical DEC solutions and known analytical solutions. The $L^2$ error is plotted against the maximum primal edge length. We note that the maximum edge length increases slightly for the non-Delaunay mesh groups due to the distortion process. The solution of all linear systems is carried out using the open source SuperLU solver \cite{superlu_ug99,superlu99}. 

\subsection{The primal 0-form Poisson equation}
\label{subsec:Primal-0-form}

The 0-form Poisson equation with Neumann boundary conditions reads
\begin{equation}
\label{eq:PoissonD2k0}
\ast \ \textrm{d} \ast \textrm{d} u = f, \ \ \ \ \ \ \ \textrm{tr} \ast \textrm{d} u = g \ \textrm{on} \ \partial \Omega.          
\end{equation} 
The above boundary condition is equivalent in vector calculus to $(\nabla u \cdot \mathbf{n} = g)$, with $\mathbf{n}$ being the outward boundary unit normal vector. Choosing to define the 0-form $u$ on the primal nodes, the discrete equations take the form  
\begin{equation}
\label{eq:PoissonD2k0-primal-DEC-0}
\left[\ast_0^{-1} [-\textrm{d}_0^T] \ast_1 \textrm{d}_0 \right] u = f - \ast_0^{-1} \textrm{d}_b g,          
\end{equation}
where the $\textrm{d}_b$ operator complements the action of the exterior derivative $[-\textrm{d}_0^T]$ for the dual cells whose boundary includes primal edges \cite{mohamed2016discrete}, and $g$ is the given Neumann boundary condition. Such a formulation is very similar to the finite element formulation of the scalar Poisson equation where Neumann boundary conditions appear as natural boundary conditions represented by a force vector.
An analytical solution that satisfies Eq. \eqref{eq:PoissonD2k0} is $u = \cos(\pi x) \cos(\pi y)$, which after substituting in Eq. \eqref{eq:PoissonD2k0} gives the right hand side as $f = -2 \pi^2 \cos(\pi x) \cos(\pi y)$. For a unit square simulation domain, the analytical solution implies zero Neumann boundary condition; i.e $g = 0$, and the final linear system to be solved is 
\begin{equation}
\label{eq:PoissonD2k0-primal-DEC-1}
\left[\ast_0^{-1} [-\textrm{d}_0^T] \ast_1 \textrm{d}_0 \right] u = f,          
\end{equation}
where $u$ is specified at one interior node in order to get a unique solution for the system.

\begin{figure}
	\centering
	\includegraphics[width=0.49\textwidth]{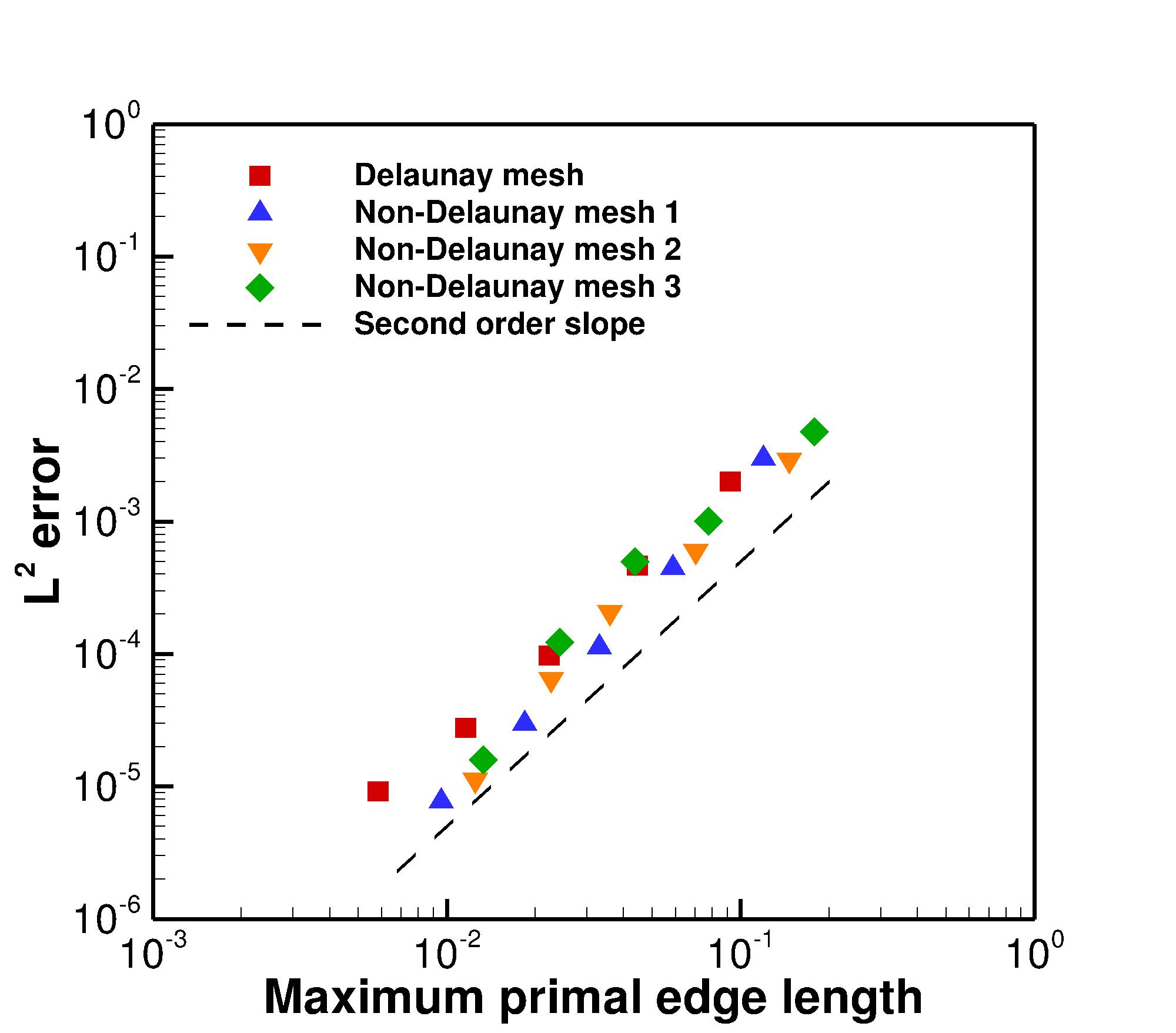}
	\caption{The $L^2$ error convergence for the primal 0-form Poisson solution with Neumann boundary conditions.}
	\label{Fig:PrimalPoissonN2k0-Convergence}
\end{figure}

Fig. \ref{Fig:PrimalPoissonN2k0-Convergence} shows the $L^2$ norm error for $u$ versus the maximum length of the primal edges for each mesh. The $L^2$ error converges in a second order fashion as expected for all mesh types. At each mesh refinement level, a slight difference in the error values exists between the Delaunay versus the non-Delaunay groups, which is attributed  to the increase in the maximum primal edges length for the non-Delaunay mesh groups. The current results indicate that the signed diagonal Hodge star operators do produce the expected correct numerical solutions even for non-Delaunay triangulations.

\subsection{The dual 0-form Poisson equation}
\label{subsec:Dual-0-form}

The 0-form Poisson equation \eqref{eq:PoissonD2k0}, solved in the previous subsection, is solved again in this subsection but with the 0-form $u$ being defined on the dual nodes. A common physical application for the dual 0-form Poisson equation is to calculate the pressure field as a post-processing step after solving the incompressible Navier-Stokes equations. The discrete dual 0-form Poisson equation takes the form  
\begin{equation}
\label{eq:PoissonD2k0-dual-0}
\left[\ast_2 \textrm{d}_1 \ast_1^{-1} \textrm{d}_1^T\right] u = f.
\end{equation}

The discrete dual 0-form Poisson equation \eqref{eq:PoissonD2k0-dual-0} with the boundary condition $\textrm{tr} \ast \textrm{d} u = g$ can also be interpreted as a primal 2-form Poisson equation by replacing $u$ by $\tilde{u} = \ast_2^{-1} u$ and $\tilde{f} = \ast_2^{-1} f$. This interpretation corresponds to the 2-form case in the boundary de Rham complex in \cite{arnold2010finite}. When the 0-form $u$ is defined on the dual nodes, the Neumann boundary condition ($\textrm{tr} \ast \textrm{d} u = g$) in Eq. \eqref{eq:PoissonD2k0} resembles a flux condition across the boundary primal edges. In order to implement such a condition, we split the exterior derivative matrix $\textrm{d}_1$ in Eq. \eqref{eq:PoissonD2k0-dual-0} as $\textrm{d}_1 = \textrm{d}_1^{\prime} + \textrm{d}_1^{\prime \prime}$, where the entries corresponding to the boundary primal edges are separated in the matrix $\textrm{d}_1^{\prime \prime}$. The discrete equations including the Neumann boundary condition $g$ are therefore
\begin{equation}
\label{eq:PoissonD2k0-dual-1}
\left[\ast_2 \textrm{d}_1^{\prime} \ast_1^{-1} \textrm{d}_1^T \right] u = f - \ast_2 \textrm{d}_1^{\prime \prime}  g.
\end{equation}   
Considering the same analytical solution $u = \cos(\pi x) \cos(\pi y)$, and a unit square simulation domain, zero Neumann boundary condition is therefore imposed; i.e, $g = 0$. The final linear system solved is 
\begin{equation}
\label{eq:PoissonD2k0-dual-2}
\left[\ast_2 \textrm{d}_1^{\prime} \ast_1^{-1} \textrm{d}_1^T \right] u = f.
\end{equation}

\begin{figure}
	\centering
	\includegraphics[width=0.49\textwidth]{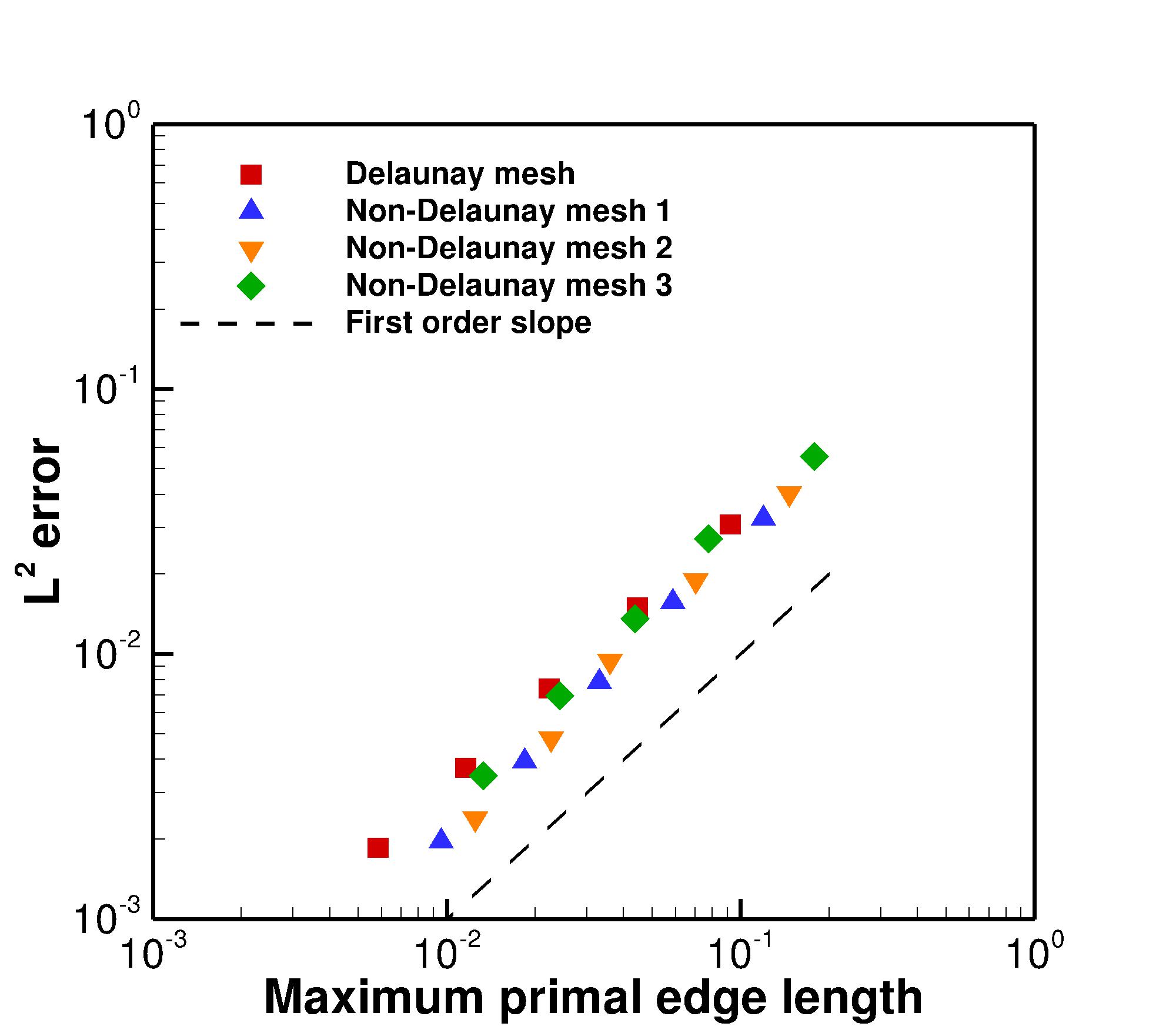}
	\caption{The $L^2$ error convergence for the dual 0-form Poisson equation with Neumann boundary conditions.}
	\label{Fig:DialPoissonD2k0-Convergence}
\end{figure}

Fig. \ref{Fig:DialPoissonD2k0-Convergence} shows the $L^2$ error convergence for the dual 0-form Poisson equation. The error converges in a first order fashion for all meshes. This first order convergence is different from the second order convergence observed before in Section \ref{subsec:Primal-0-form} for the primal 0-form Poisson equation. During the current dual 0-form formulation, the 0-forms $u$ defined on the dual nodes are numerically considered to be constant over the triangles. Accordingly, only a first order convergence is expected. This is also consistent with the first order convergence in lowest order FEEC solution of the 2-form Poisson equation \cite{arnold2010finite}. Insofar as the error is concerned, the slight increase in the non-Delaunay meshes error seems to be consistent with the corresponding increase in the maximum primal edges length. Correct numerical solutions are therefore obtained regardless of the non-Delaunay triangles ratio while using the circumcentric dual mesh and the signed diagonal Hodge operators.      

\subsection{The 1-form Poisson equation}
\label{subsec:1-form-Poisson}

A 1-form Poisson equation with Dirichlet boundary conditions is expressed in exterior calculus notation as
\begin{equation}
\label{eq:PoissonD2k1}
\textrm{d} \ast \textrm{d} \ast u + \ast \ \textrm{d} \ast \textrm{d}  u = f, \ \ \ \ \ \ \ \textrm{tr} \ u = g, \ \textrm{tr} \ast u = h \ \ \textrm{on} \ \partial \Omega.  
\end{equation}
Choosing to define the 1-form $u$ on the primal edges, the DEC formulation takes the form
\begin{equation}
\label{eq:PoissonD2k1-DEC}
\left[\textrm{d}_0 \ast_0^{-1} [-\textrm{d}_0^T] \ast_1 + \ast_1^{-1} \textrm{d}_1^T \ast_2 \textrm{d}_1 \right] u =  f - \textrm{d}_0 \ast_0^{-1} \textrm{d}_b h.
\end{equation}
The boundary condition $\textrm{tr} \ast u = h$ appears as a force vector on the right hand side similar to Eq. \eqref{eq:PoissonD2k0-primal-DEC-0}. The Dirichlet boundary condition $\textrm{tr} \ u = g$ applied on each boundary edge is implemented by multiplying the corresponding column in the stiffness matrix by the imposed boundary condition value, subtracting the result from the right hand side force vector, and finally removing both the corresponding row and column from the mass matrix. We pose the problem on a simply connected domain in order to avoid the additional requirement of computing harmonic forms. An analytical solution for the 1-form Poisson equation on a unit square domain is 
\begin{equation}
\label{eq:PoissonD2k1-analy}
u = x (1-x) y^2 (1-y)^2 dx + y (1-y) x^2 (1-x)^2 dy.  
\end{equation}
Such an analytical solution results in zero values for both boundary conditions.

\begin{figure}
	\centering
	\includegraphics[width=0.49\textwidth]{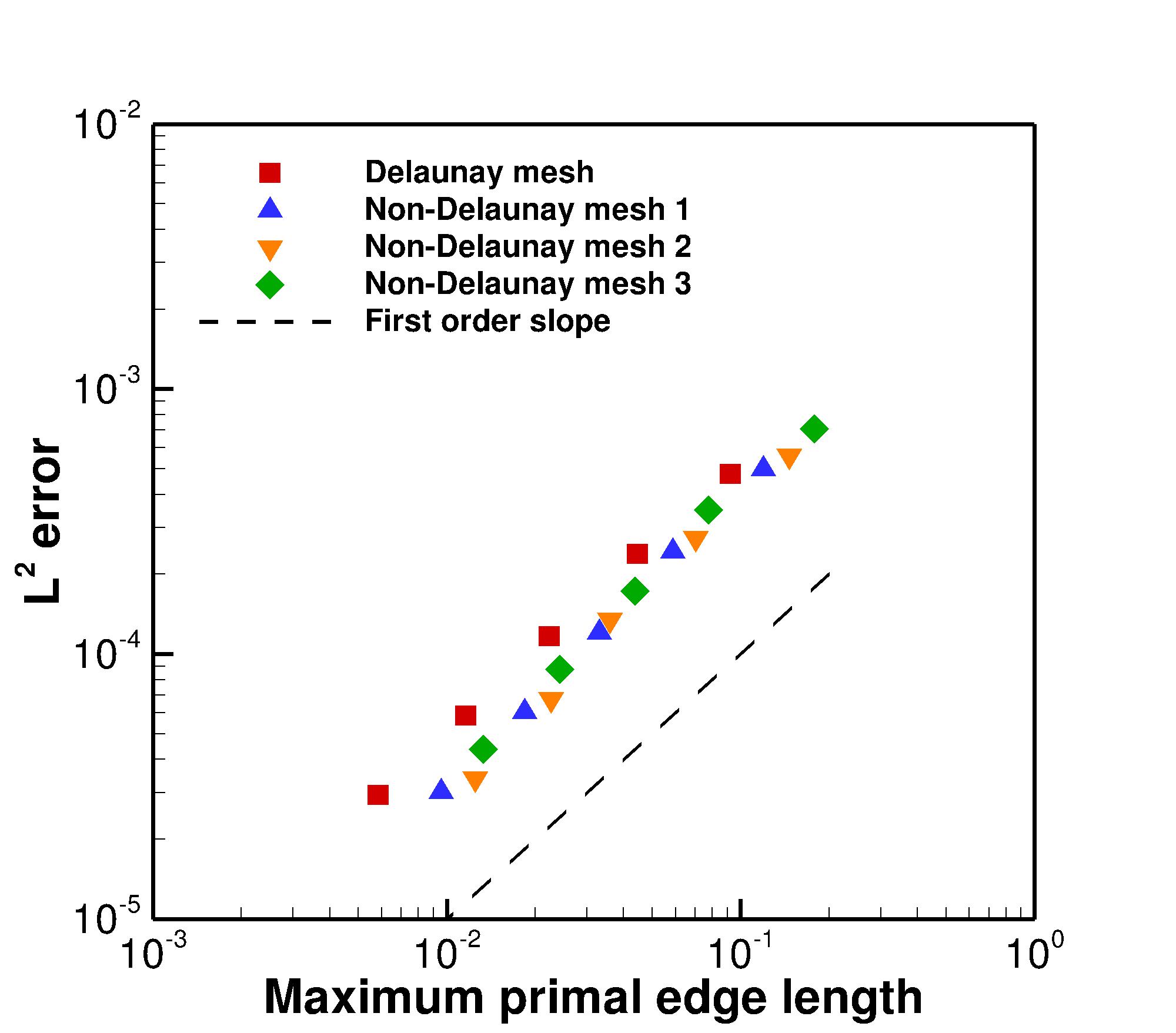}
	\caption{The $L^2$ error convergence for the interpolated 1-form $u$ after solving the 1-form Poisson equation.}
	\label{Fig:PoissonD2k1-Convergence}
\end{figure}

The $L^2$ error convergence for the 1-form $u$ is shown in Fig. \ref{Fig:PoissonD2k1-Convergence}. For all test cases, the error values converge with a first order rate as expected. Again, the error values for all meshes vary almost consistently with the maximum primal edges length. The results further indicate the ability of the signed diagonal Hodge operators to provide correct DEC solutions even for meshes having a significant fraction of non-Delaunay triangles.     

\subsection{The incompressible Navier-Stokes equations}
\label{subsec:Incomp-Navier-Stokes}

The DEC formulation of the incompressible Navier-Stokes equations is briefly presented in this subsection, and the reader may refer to \cite{mohamed2016discrete} for more details. For a homogeneous fluid having a unit density and no body forces, the incompressible flow of the fluid over a flat surface is governed by     
\begin{subequations}
	\label{eq:NS_DG-00}
	\begin{align}
	\label{eq:NS_DG-00-a}
	&\frac{\partial \mathbf{u}^\flat}{\partial t} - \nu \ast \textrm{d} \ast \textrm{d} \mathbf{u}^\flat + \ast(\mathbf{u}^\flat \wedge \ast \textrm{d} \mathbf{u}^\flat) + \textrm{d}p^d =0,   \\
	\label{eq:NS_DG-00-b}
	&\ast \textrm{d} \ast \mathbf{u}^\flat = 0,
	\end{align}
\end{subequations}
with $\mathbf{u}^\flat$ being the velocity 1-form, the dynamic pressure is $p^d = p + \frac{1}{2} (\mathbf{u}^\flat(\mathbf{u}))$, and $\nu$ is the kinematic viscosity. The flow over a curved surface requires a curvature term (only in the viscous case) to be added to the momentum equation \cite{EbMa1970}. Since the presented experiments are either viscous flows over flat surfaces or inviscid flows over curved surfaces, the curvature term is not considered during the current discretization. The DEC discrete equations are obtained by first taking the exterior derivative of the momentum equation \eqref{eq:NS_DG-00-a}, substituting with the DEC discrete operators, and replacing the velocity dual 1-form $u$ by $u = \ast_1 \textrm{d}_0 \psi$, where $\psi$ is the stream function primal 0-form. The discrete system of equations to be solved is then 
\begin{multline}
\label{eq:discreteNS02}
\frac{1}{\Delta t} [-\textrm{d}^T_0] \ast_1 \textrm{d}_0 \psi^{n+1} - \nu [-\textrm{d}^T_0] \ast_1 \textrm{d}_0 \ast^{-1}_0 [-\textrm{d}^T_0] \ast_1 \textrm{d}_0 \psi^{n+1} \\ + [-\textrm{d}^T_0] \ast_1 W_v^n \ast^{-1}_0 [-\textrm{d}^T_0] \ast_1 \textrm{d}_0 \psi^{n+1}= F,
\end{multline}
where $F = \frac{1}{\Delta t} [-\textrm{d}^T_0] u^n + \nu [-\textrm{d}^T_0] \ast_1 \textrm{d}_0 \ast^{-1}_0  \textrm{d}_b v^n -  [-\textrm{d}^T_0] \ast_1 W_v^n \ast^{-1}_0  \textrm{d}_b v^n $. The matrix $W_v$, containing the tangential velocity primal 1-form $v$ values, represents the discrete wedge product action. The degrees of freedom in Eq. \eqref{eq:discreteNS02} are the stream function 0-forms $\psi$ defined on the primal mesh nodes.

A test case for incompressible Navier-Stokes equations that has a known analytical solution is the Poiseuille flow. For a rectangular simulation domain, solid wall boundary conditions are applied on the top and bottom boundaries, while parabolic in/out flow conditions are imposed on the left/right boundaries. For a unit square domain, the exact solution for the velocity vector field is $\mathbf{u} = \left(y(1-y),0\right)$, with $\nu = 1.0$. As explained earlier in \cite{mohamed2016discrete}, the velocity 1-form $L^2$-error is calculated as $\| u^{exact} - u\| = \left[ \sum_{\sigma^1} A_s \left(\frac{u^{exact} - u}{|\sigma^1|}\right)^2 \right]^{1/2}$, where $A_s$ is the primal edge $\sigma^1$ support area represented by the polygonal area covering both the primal edge and its dual. For flipped dual edges, $A_s$ is considered with a negative sign.

\begin{figure}
	\centering
	\includegraphics[width=0.49\textwidth]{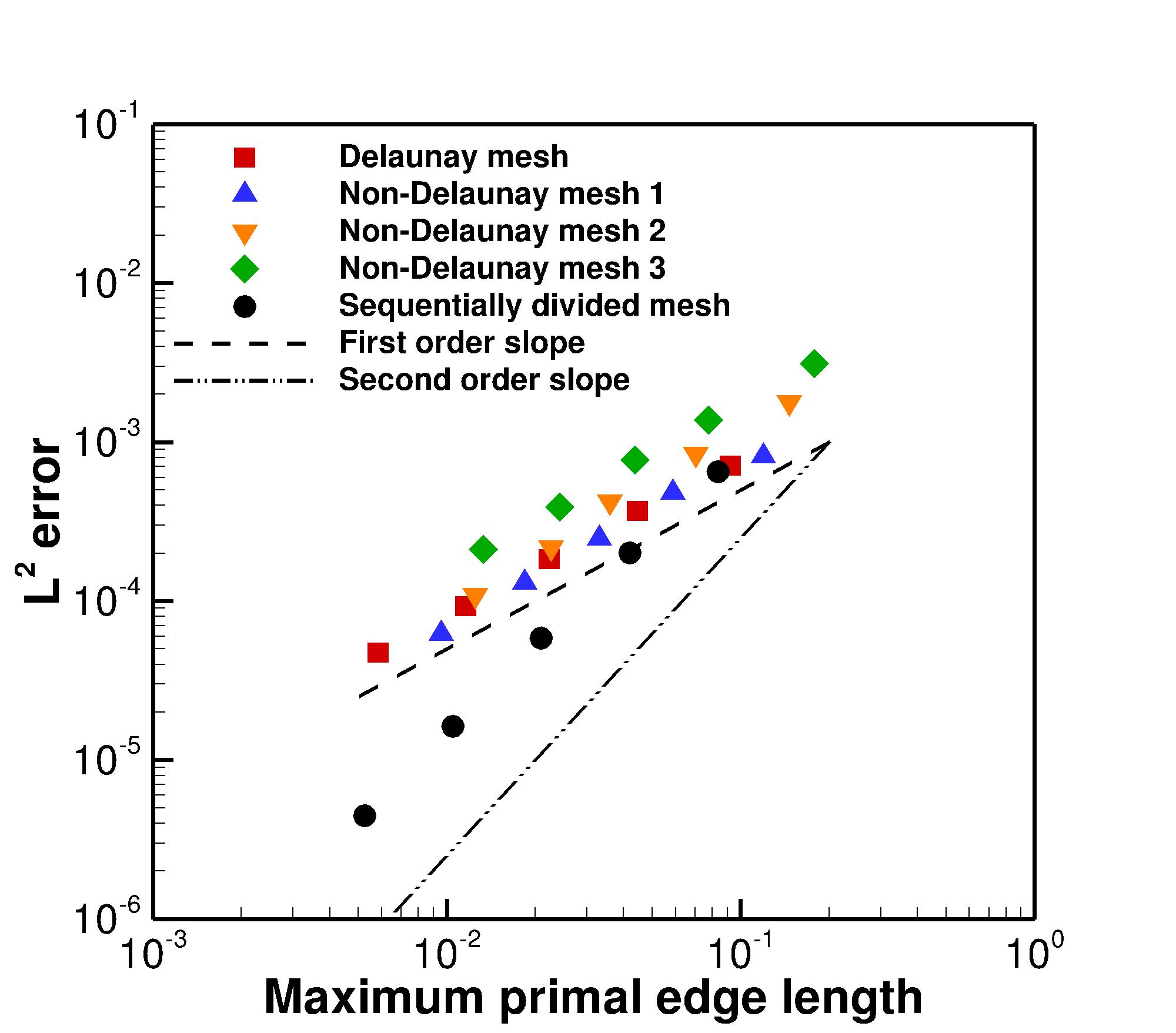}
	\caption{The velocity 1-form $L^2$ error convergence for the Poiseuille flow test case.}
	\label{Fig:Poiseuille-Convergence}
\end{figure}

Fig. \ref{Fig:Poiseuille-Convergence} shows the $L^2$ error convergence for the Poiseuille flow test case. For both Delaunay and non-Delaunay mesh groups, the error converges in a first order fashion as expected for unstructured meshes. In regards to the non-Delaunay mesh 3 group, the error increases slightly above that expected based on the maximum primal edges length. However, it does appear to converge in the expected first order fashion.

An additional group of sequentially-divided non-Delaunay meshes is considered for the incompressible Navier-Stokes solution. This group of meshes is generated by sequentially subdividing an initial coarse Delaunay mesh by adding new nodes at the edges midpoints. The percentage ratio of non-Delaunay triangles increases gradually at each refinement level and is $\sim 12 \%$ for the finest mesh. Such a group of meshes was tested earlier in \cite{mohamed2016baryHodge}, but with a barycentric dual mesh defined over the whole triangulation, and exhibited a convergence rate very close to a second order. This was attributed to the subdivision process that introduces new interior edges having their midpoints coincident with their barycentric dual edges midpoints. Such a condition was shown earlier by Nicolaides \cite{nicolaides:1992direct} to result in super-convergence. The results displayed in Fig. \ref{Fig:Poiseuille-Convergence} show that such super-convergence is retained even when the solution is carried out using a circumcentric dual mesh. Similar to the barycentric dual situation, the newly introduced interior primal edges created during the subdivision process have their midpoints coincident with the midpoints of their circumcentric dual edges. Such a property retains the super-convergence effect when using the circumcentric dual mesh, although the primal triangulation is non-Delaunay.

\subsection{Condition number analysis}
\label{subsec:Condition-number}

The previous four subsections provided numerical evidence that DEC solutions using the circumcentric dual and the diagonal Hodge star operators do provide correct results even with non-Delaunay triangulations. This subsection provides additional insight into how the stiffness matrices condition number can change when using non-Delaunay meshes. The changes in the matrices condition number can be related to the changes in the minimum volumes of the primal edges, dual edges, primal triangles and dual cells. This is expected since the Hodge star operators perform divisions by such volumes. For entities having signed volumes; e.g. the dual cells and the dual edges, whenever the minimum volume is mentioned here it is meant to be that having the smallest absolute value. 

\begin{figure}
	\centering
	\begin{subfigure}[b]{0.48\textwidth}
		\centering
		\includegraphics[width=\textwidth]{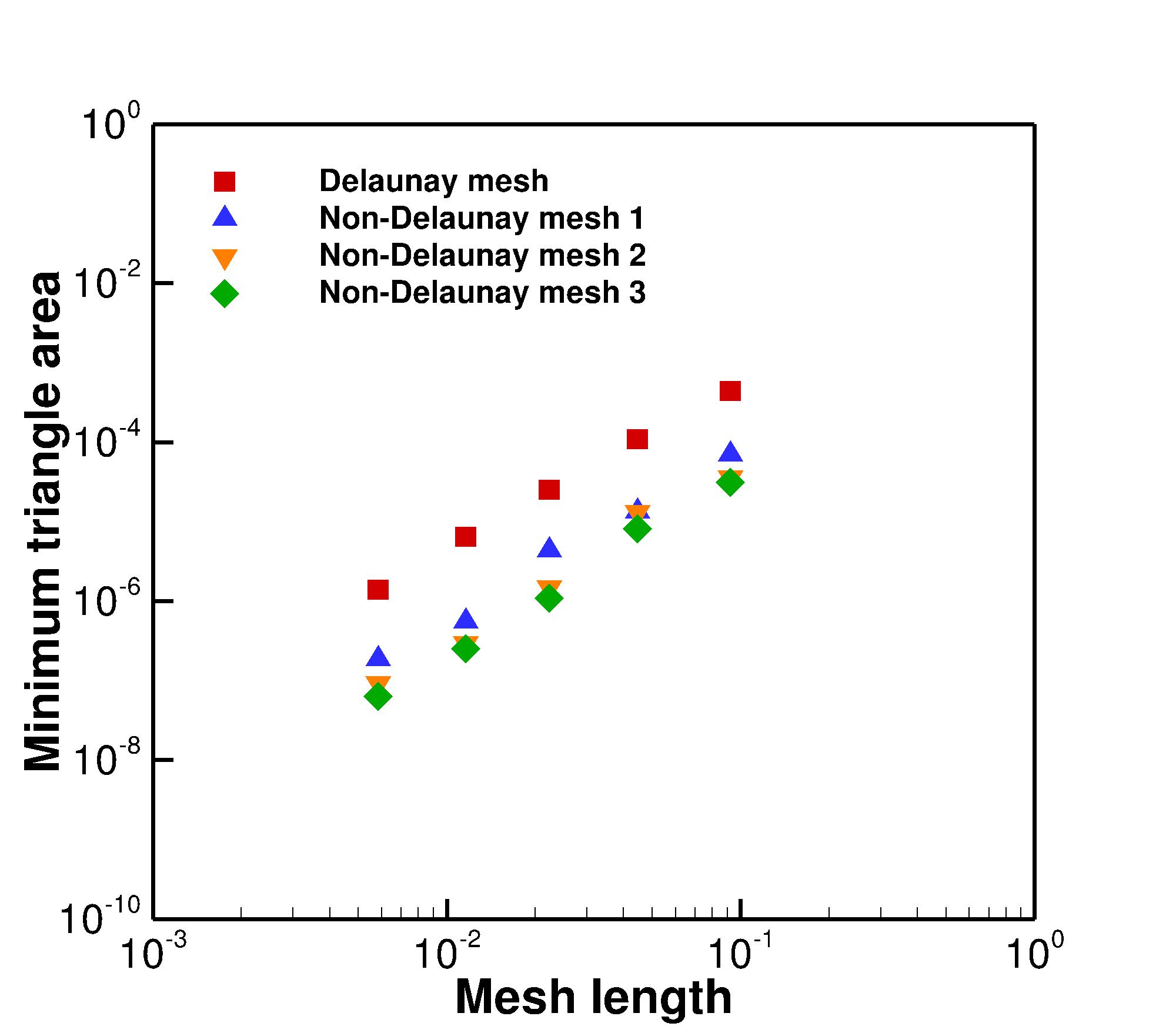}
		\caption{}
		\label{subFig:areas-a}
	\end{subfigure} %
	\begin{subfigure}[b]{0.48\textwidth}
		\centering
		\includegraphics[width=\textwidth]{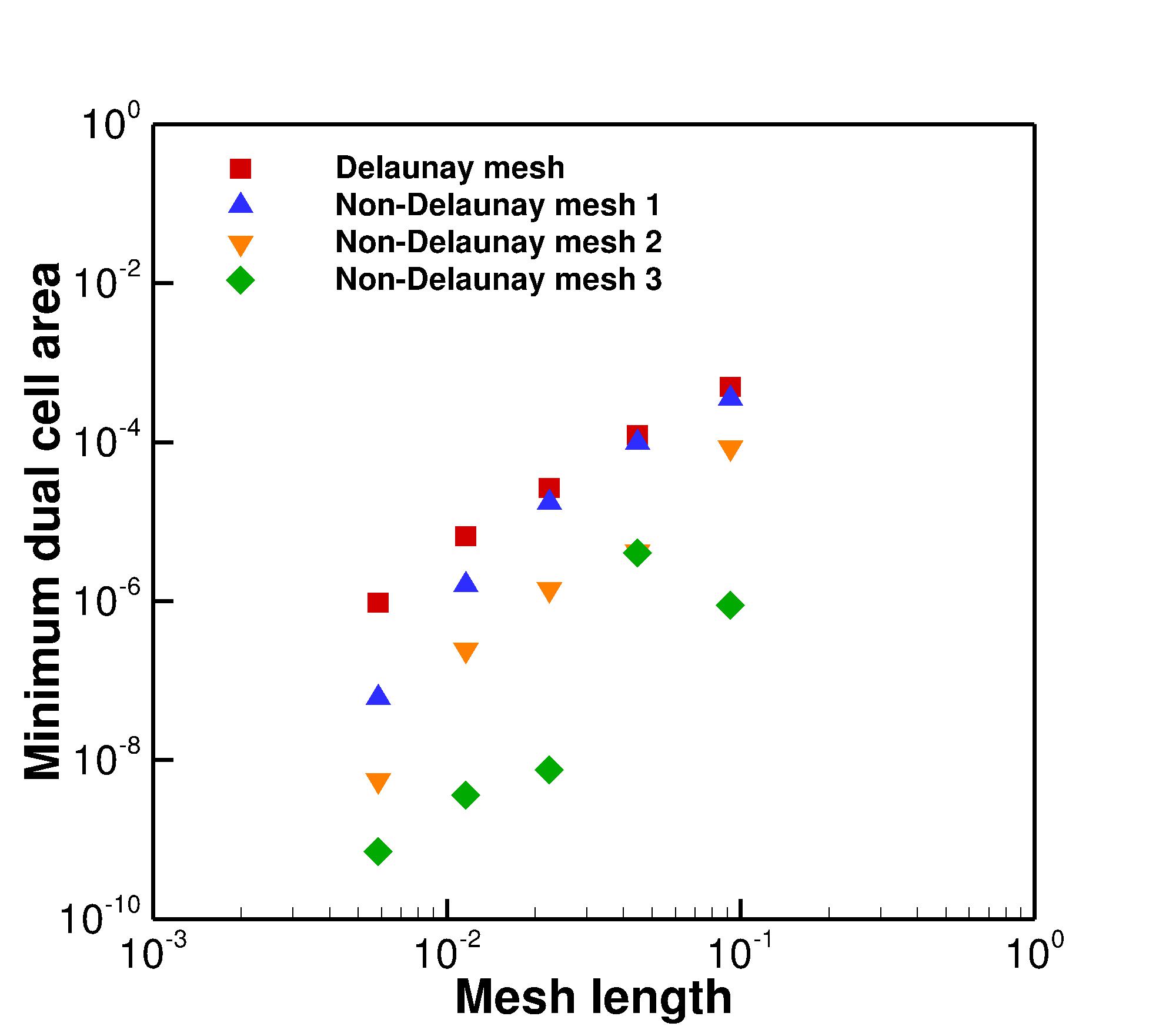} 
		\caption{}
		\label{subFig:areas-b}
	\end{subfigure} %
	\caption{(\subref{subFig:areas-a}) The minimum primal triangle area for all considered meshes. (\subref{subFig:areas-b}) The minimum (as an absolute value) dual cell area for all meshes. For all data points, the x-axis values are the maximum primal edge length for the Delaunay mesh group in order to avoid any unnecessary data shift in the x-axis. The data points for the sequentially divided non-Delaunay mesh group are almost coincident with the Delaunay mesh group data points in both figures, and were therefore excluded for better visibility.}
	\label{Fig:areas}
\end{figure}

Figure \ref{Fig:areas} shows the minimum triangle area and minimum dual cell area for all considered meshes. In comparison to the Delaunay mesh group, the non-Delaunay meshes exhibit a decrease in the minimum triangle area. This is due to the distortion process which results in triangles having relatively large aspect ratio, and therefore a relatively smaller area. The decrease in the minimum triangle area is almost by one order of magnitude. The minimum dual cell areas, on the other hand, exhibit a more significant decrease of nearly four orders of magnitude for the non-Delaunay triangles. This is attributed to the negative-volume area sectors that seem to significantly decrease the dual cells net area. The Hodge star operator $\ast_0^{-1}$ includes a division by the dual cells areas. Accordingly, the stiffness matrices including $\ast_0^{-1}$ are expected to have higher condition number for the non-Delaunay mesh groups. The changes in the minimum primal and dual edge length will be addressed during the condition number results discussion. It will be shown that while the minimum primal edge length decreases by almost one order of magnitude for the non-Delaunay meshes, the minimum dual edge length does not change significantly. This is because tiny dual edges can and do exist even for Delaunay meshes. Considering the sequentially divided non-Delaunay mesh group, the values of both the minimum triangle area and the minimum dual cell area are almost coincident with that for the Delaunay mesh group, and were therefore excluded from Fig. \ref{Fig:areas} for better visibility. 

\begin{figure}
	\centering
	\begin{subfigure}[b]{0.32\textwidth}
		\centering
		\includegraphics[width=\textwidth]{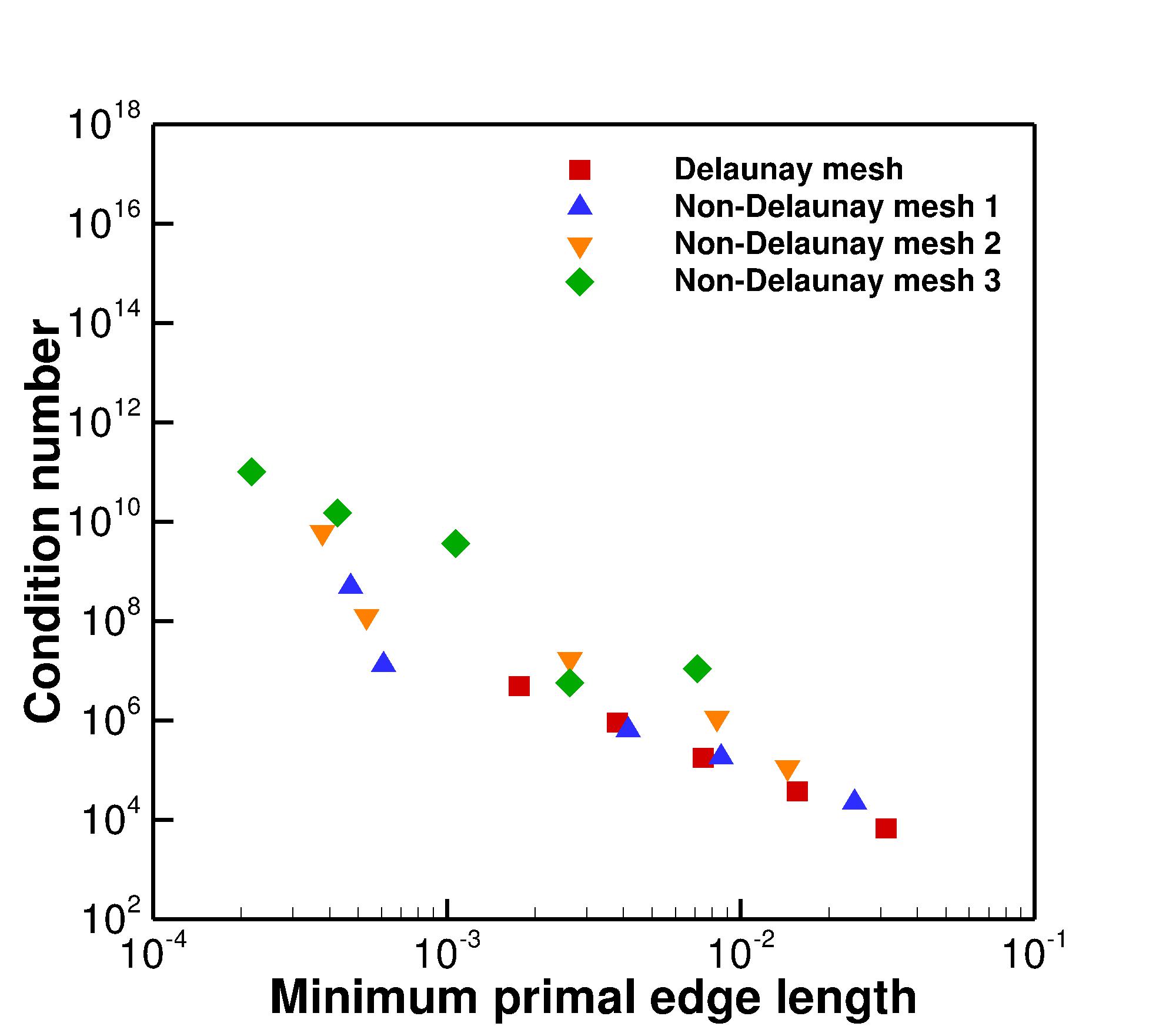}
		\caption{}
		\label{subFig:condNum-a}
	\end{subfigure} %
	\begin{subfigure}[b]{0.32\textwidth}
		\centering
		\includegraphics[width=\textwidth]{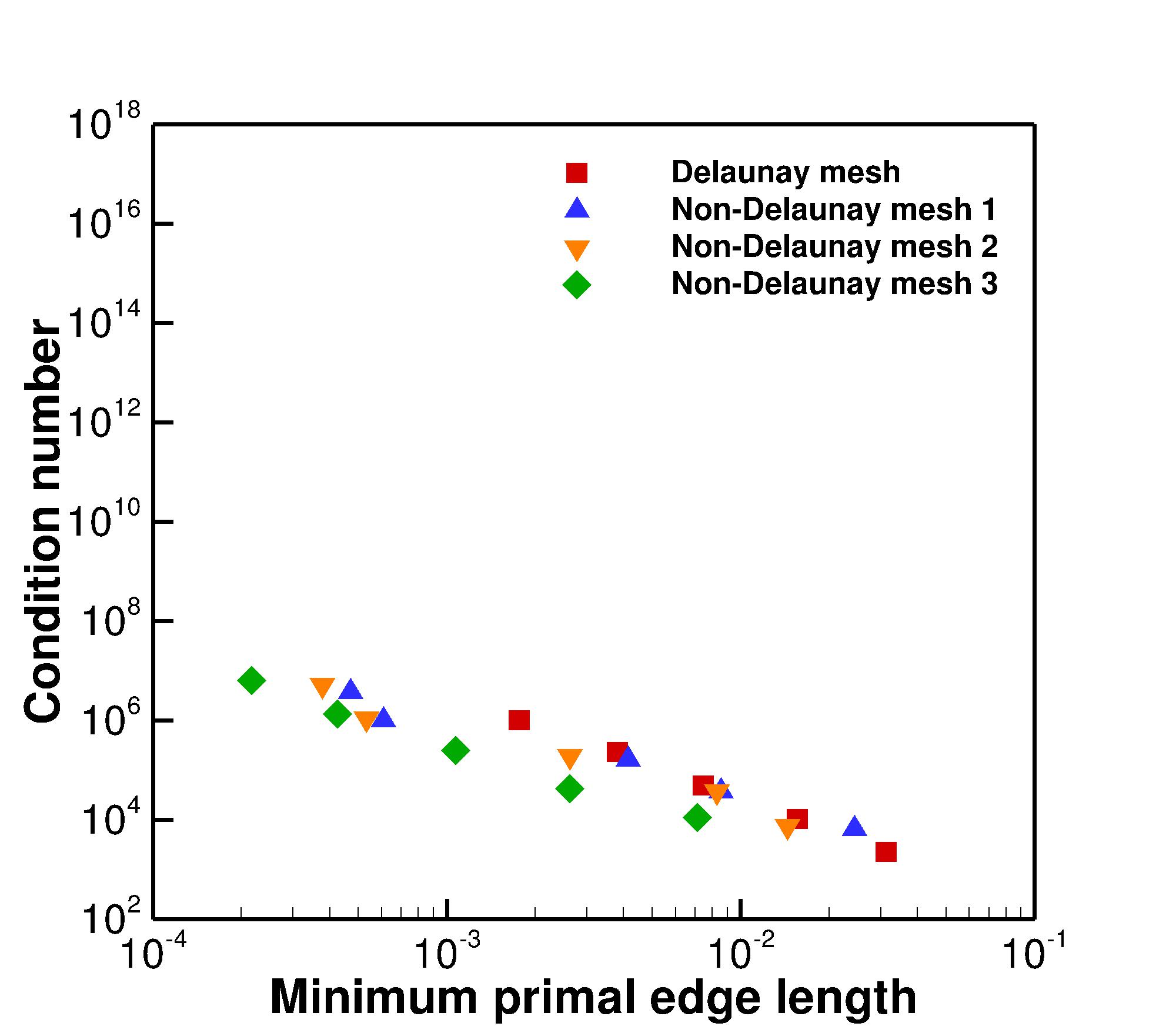} 
		\caption{}
		\label{subFig:condNum-b}
	\end{subfigure} %
	\begin{subfigure}[b]{0.32\textwidth}
		\centering
		\includegraphics[width=\textwidth]{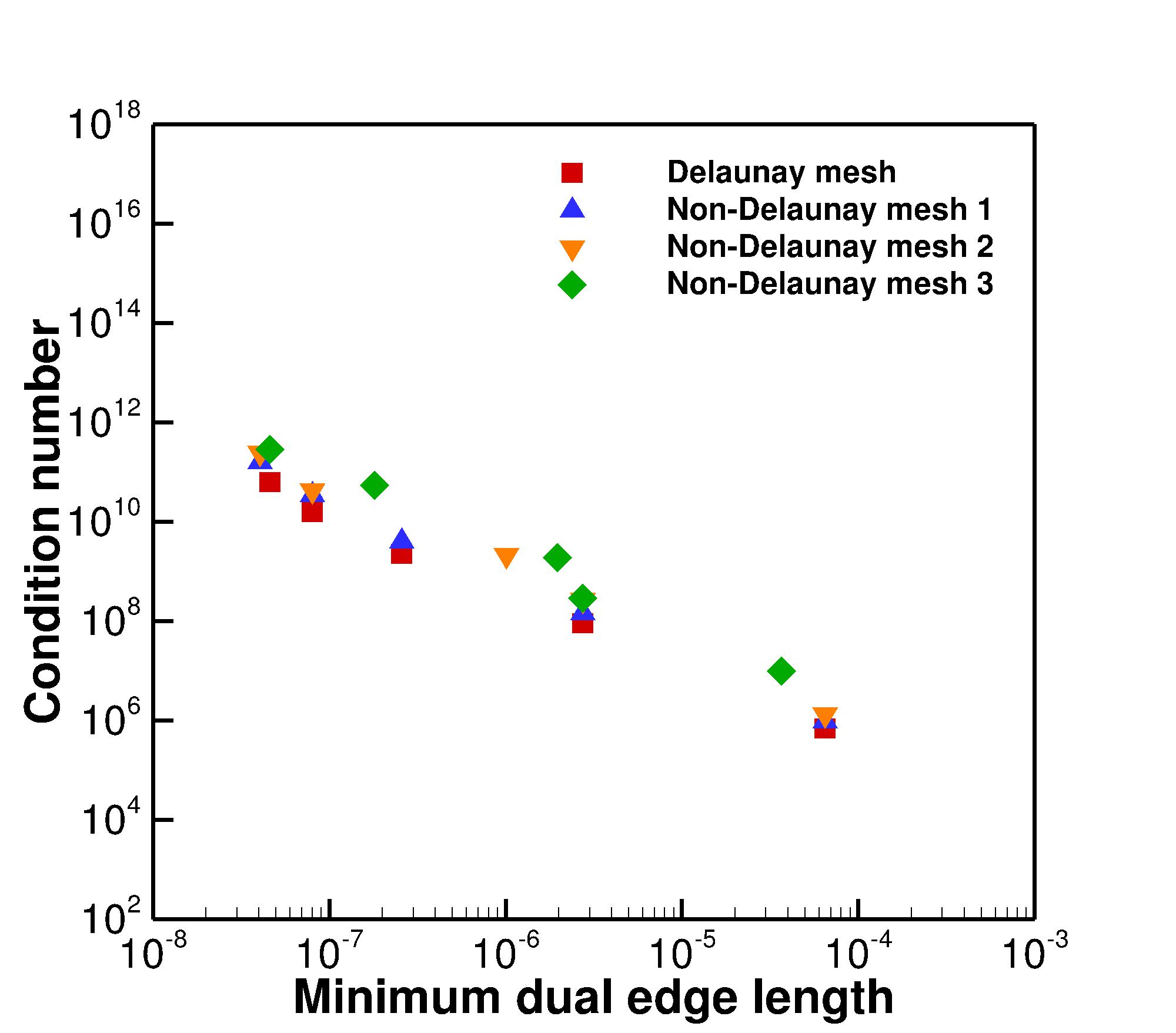} 
		\caption{}
		\label{subFig:condNum-c}
	\end{subfigure} %
	\begin{subfigure}[b]{0.32\textwidth}
		\centering
		\includegraphics[width=\textwidth]{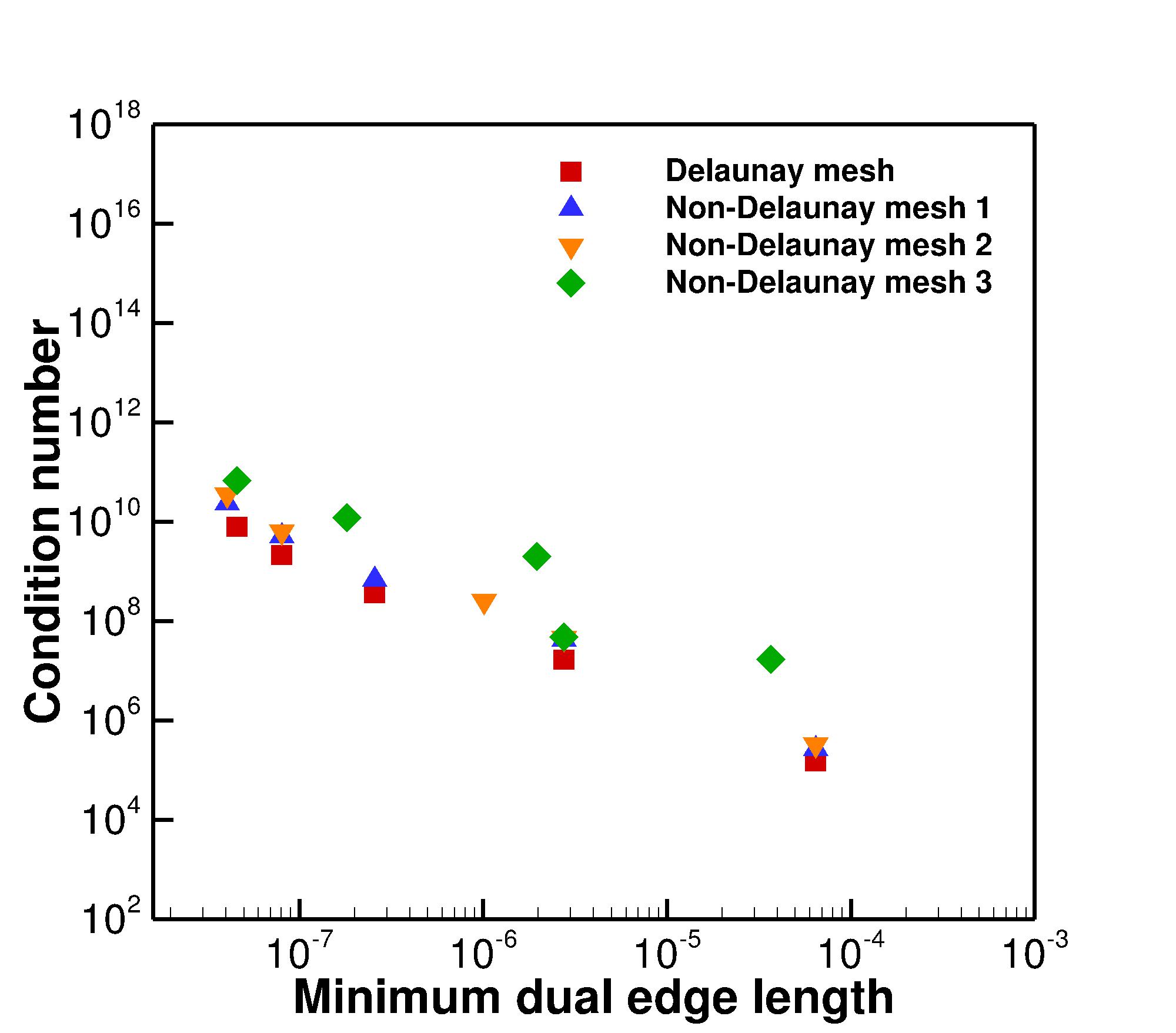} 
		\caption{}
		\label{subFig:condNum-d}
	\end{subfigure} %
	\begin{subfigure}[b]{0.32\textwidth}
		\centering
		\includegraphics[width=\textwidth]{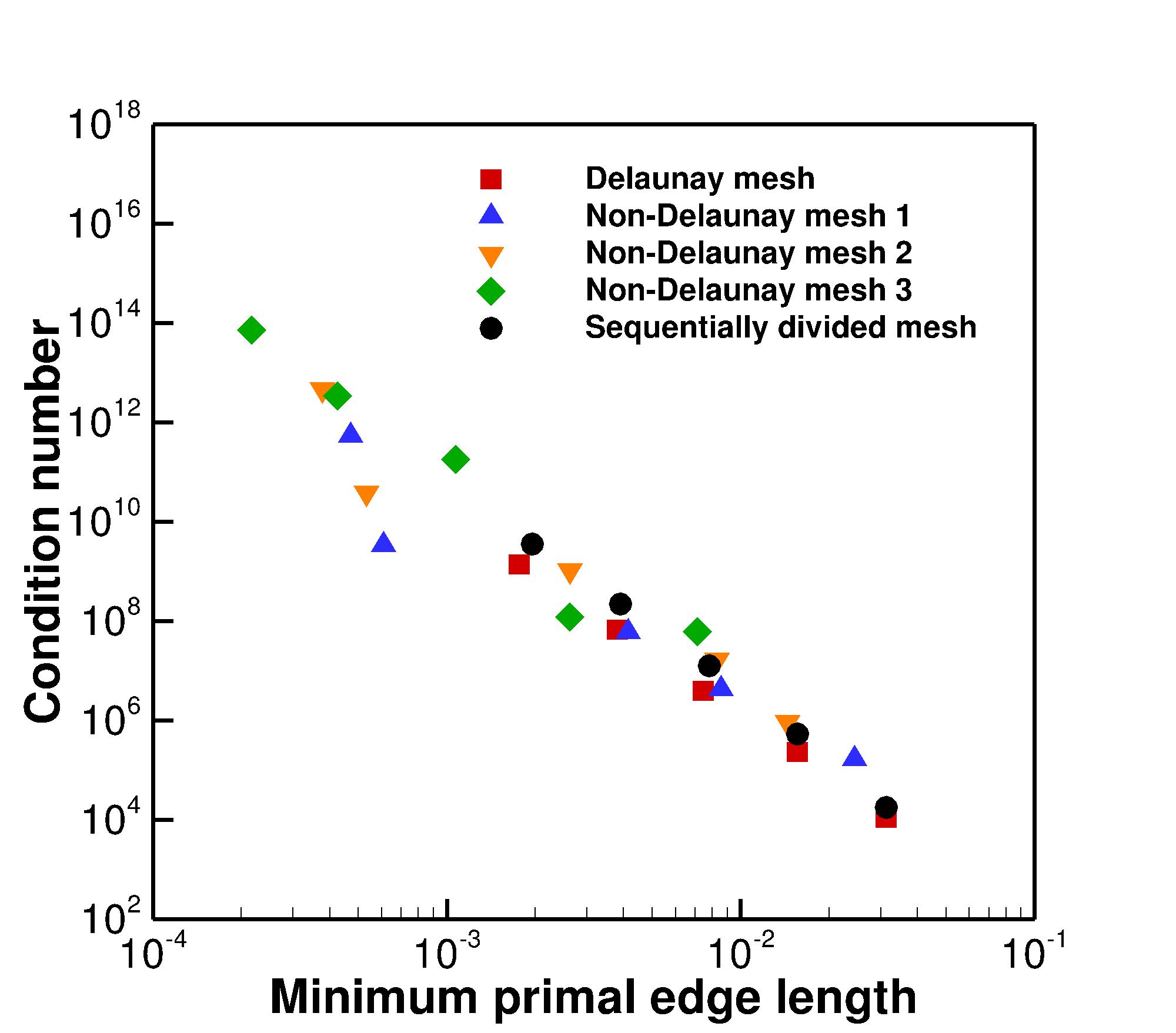} 
		\caption{}
		\label{subFig:condNum-e}
	\end{subfigure} %
	\caption{The condition number for: (\subref{subFig:condNum-a}) The matrix $\left[\ast_0^{-1} [-\textrm{d}_0^T] \ast_1 \textrm{d}_0 \right]$ used in the primal 0-form Poisson equation solved in Section \ref{subsec:Primal-0-form}.  (\subref{subFig:condNum-b}) The matrix $\left[ [-\textrm{d}_0^T] \ast_1 \textrm{d}_0 \right]$ used in the primal 0-form Poisson equation after being multiplied by $\ast_0$. (\subref{subFig:condNum-c}) The matrix $\left[\ast_2 \textrm{d}_1^{\prime} \ast_1^{-1} \textrm{d}_1^T \right]$ used in the dual 0-form Poisson equation solved in Section \ref{subsec:Dual-0-form}. (\subref{subFig:condNum-d}) The matrix $\left[\textrm{d}_0 \ast_0^{-1} [-\textrm{d}_0^T] \ast_1 + \ast_1^{-1} \textrm{d}_1^T \ast_2 \textrm{d}_1 \right]$ used in the primal 1-form Poisson equation solved in Section \ref{subsec:1-form-Poisson}. (\subref{subFig:condNum-e}) The incompressible Navier Stokes stiffness matrix solved in Section \ref{subsec:Incomp-Navier-Stokes}.}
	\label{Fig:condNum}
\end{figure}

The stiffness matrices condition numbers for the four solved partial differential equations are shown in Fig. \ref{Fig:condNum}. First, Fig. \ref{subFig:condNum-a} shows the condition number of the $\left[\ast_0^{-1} [-\textrm{d}_0^T] \ast_1 \textrm{d}_0 \right]$ matrix which appears in the primal 0-form Poisson equation solved in Section \ref{subsec:Primal-0-form}. The condition number is plotted versus the minimum primal edge length because the Hodge operator $\ast_1$ involves division by the primal edges length. The minimum primal edge length decreases for the non-Delaunay meshes, due to the distortion process, and therefore the matrix condition number increases. At the same refinement level, the condition number for non-Delaunay meshes is almost four orders of magnitude larger than that for the Delaunay meshes. Another factor, however, that can impact the condition number is the division by the dual cells areas implied in the Hodge star operator $\ast_0^{-1}$. This can be verified by multiplying Eq. \eqref{eq:PoissonD2k0-primal-DEC-1} by $\ast_0$, which removes $\ast_0^{-1}$ from the left hand side. Solving the equation $[-\textrm{d}_0^T] \ast_1 \textrm{d}_0 u = \ast_0 f$ results in error values identical to these plotted in Fig. \ref{Fig:PrimalPoissonN2k0-Convergence} but with a much lower condition number for the stiffness matrix as shown in Fig. \ref{subFig:condNum-b}. This indicates how detrimental to the condition number the division by the dual cells areas can be for non-Delaunay triangulations. However, it is worth emphasizing that in spite of the relatively high condition number for the stiffness matrix in the original Eq. \eqref{eq:PoissonD2k0-primal-DEC-1}, the linear system was successfully solved using the open source SuperLU solver. The multiplication by $\ast_0$ can be considered as a matrix preconditioning, which can be beneficial when considering iterative linear solvers.                   

The stiffness matrix condition number for the dual 0-form and the primal 1-form Poisson equations are shown in Figs \ref{subFig:condNum-c} and \ref{subFig:condNum-d}, respectively. The condition numbers in both figures are plotted against the minimum dual edge length since the division by it is implied in the Hodge star operator $\ast_1^{-1}$. The primal 1-form Poisson equation includes division by both the primal edges length (in the $\ast_1$ operator) and the dual edges length (in the $\ast_1^{-1}$ operator). However, the minimum dual edge length values are significantly smaller than the minimum primal edge length values, as can be observed in the x-axis numbers in Figs. \ref{subFig:condNum-a} and \ref{subFig:condNum-d}. This suggests that division by the minimum dual edges length can be a dominating factor, and therefore the condition numbers in Fig. \ref{subFig:condNum-d} are plotted against the minimum dual edge length. Both Figs \ref{subFig:condNum-c} and \ref{subFig:condNum-d} reveal an insignificant change in the minimum dual edges length for the non-Delaunay meshes, in comparison to the Delaunay meshes. This is expected since tiny dual edges can exist even for high quality Delaunay triangulations. The stiffness matrix condition numbers for both the dual 0-form and the primal 1-form Poisson equations do not exhibit significant changes for the non-Delaunay versus the Delaunay mesh groups, at a given refinement level. The condition numbers seem to mostly be governed by the minimum dual edges length.   

Finally, Fig. \ref{subFig:condNum-e} shows the condition number for the incompressible Navier-Stokes stiffness matrix in Eq. \eqref{eq:discreteNS02}. At the same refinement level, the condition number for the non-Delaunay meshes experiences a significant increase in comparison to the Delaunay meshes. This can be caused by the decrease in both the minimum primal edge length and the minimum dual cell area through the $\ast_1$ and $\ast_0^{-1}$ operators, respectively. Equation \eqref{eq:discreteNS02} includes two $\ast_1$ operators, which suggests that the division by the primal edges length plays a dominant role in the matrix condition number increase. Even with condition number values close to $10^{14}$, the SuperLU solver could correctly solve all linear systems,  and hence led to the convergent solutions shown in Fig. \ref{Fig:Poiseuille-Convergence}. For the sequentially divided non-Delaunay mesh group, since the values of both the minimum dual cells areas and the minimum primal edges length are very close to these of the Delaunay mesh group, the condition numbers changed very little when compared with the Delaunay mesh group, as shown in Fig. \ref{subFig:condNum-e}.

The aforementioned results suggest that the two main factors driving the increase in the stiffness matrix condition numbers for non-Delaunay triangulations can be the decrease in both the minimum primal edges length and the dual cells area. However, these two factors are not characteristic features for all non-Delaunay triangulations, in comparison to Delaunay triangulations. It is evident from Fig. \ref{Fig:mesh} that the negative-volume dual cell sector can be relatively small in comparison to the positive-volume sectors, but its effect will increase as the two non-Delaunay triangles are further squeezed. Therefore, the observed significant decrease in the minimum dual cells areas can be mostly due to the distortion process we carried out aiming deliberately to generate highly distorted non-Delaunay triangulations having relatively high aspect ratios. The decrease in the minimum primal edges lengths can also be a consequence of the distortion process, and is not a characteristic feature for non-Delaunay triangulations. The non-Delaunay meshes considered in this research were meant to serve as relatively extreme cases of non-Delaunay triangulations. Such meshes may not be encountered in practical applications where the generated meshes are not deliberately intended to be of inferior quality. In such practical cases, the stiffness matrices condition number for non-Delaunay triangulations can be not significantly different from that expected for Delaunay triangulations.

\subsection{Flow simulations on non-Delaunay curved surfaces}
\label{subsec:curved-surfaces}

The non-Delaunay triangulations are likely to be encountered in mainly two situations. First, non-Delaunay triangulations may occur whenever a subdivision is carried out for a given triangulation. Even when the initial triangulation is Delaunay, the existence of any obtuse-angled triangles results in a non-Delaunay mesh after the subdivision process. An example of such a mesh is the sequentially subdivided non-Delaunay mesh group that was used in Section \ref{subsec:Incomp-Navier-Stokes}. In terms of the numerical error, the sequentially subdivided mesh group exhibited super-convergence during the solution of the incompressible Navier Stokes equation, as was shown in Fig. \ref{Fig:Poiseuille-Convergence}. In addition, this mesh group did not exhibit almost any increase in the stiffness matrix condition number, in comparison to the Delaunay mesh group. These two observations suggest that although being non-Delaunay, DEC solutions using the circumcentric dual can perform very efficiently with such class of subdivided meshes, and potentially achieve super-converging solutions for some problems.         

The second situation in which non-Delaunay triangulations may be encountered is on curved surfaces. Although various open source libraries are capable of generating Delaunay triangulations for both 2D flat surfaces and 3D domains, generating Delaunay meshes on curved surfaces remains a complicated process. In this subsection, we present a test case where a non-Delaunay triangulation is generated on a curved surface through a very simple scheme, and the behavior of DEC solution over this surface is investigated.

\begin{figure}
	\centering
	\begin{subfigure}[b]{0.32\textwidth}
		\centering
		\includegraphics[width=\textwidth]{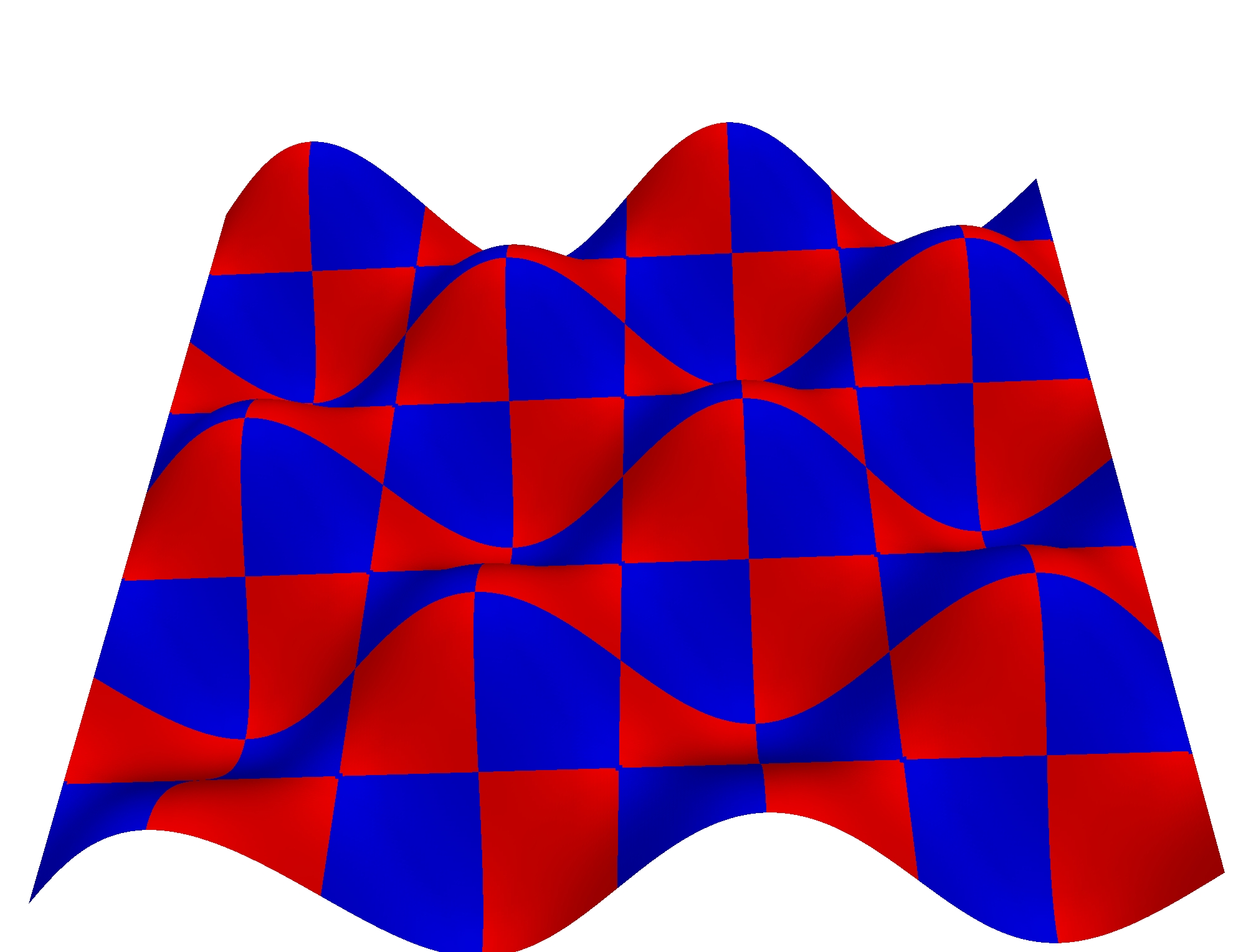}
		\caption{}
		\label{subFig:doubleShearCurved-a}
	\end{subfigure} %
	\begin{subfigure}[b]{0.32\textwidth}
		\centering
		\includegraphics[width=\textwidth]{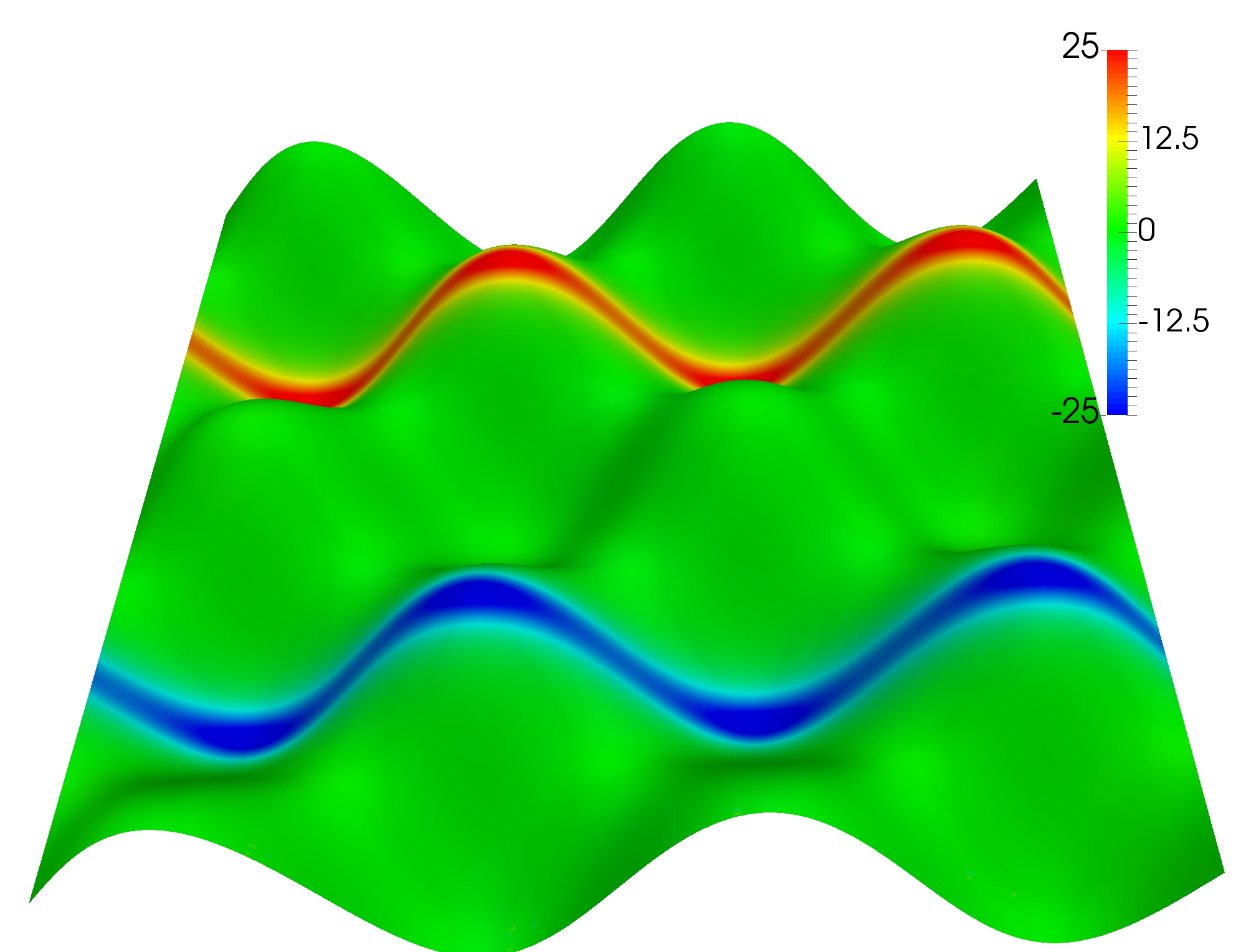} 
		\caption{}
		\label{subFig:doubleShearCurved-b}
	\end{subfigure} %
	\begin{subfigure}[b]{0.32\textwidth}
		\centering
		\includegraphics[width=\textwidth]{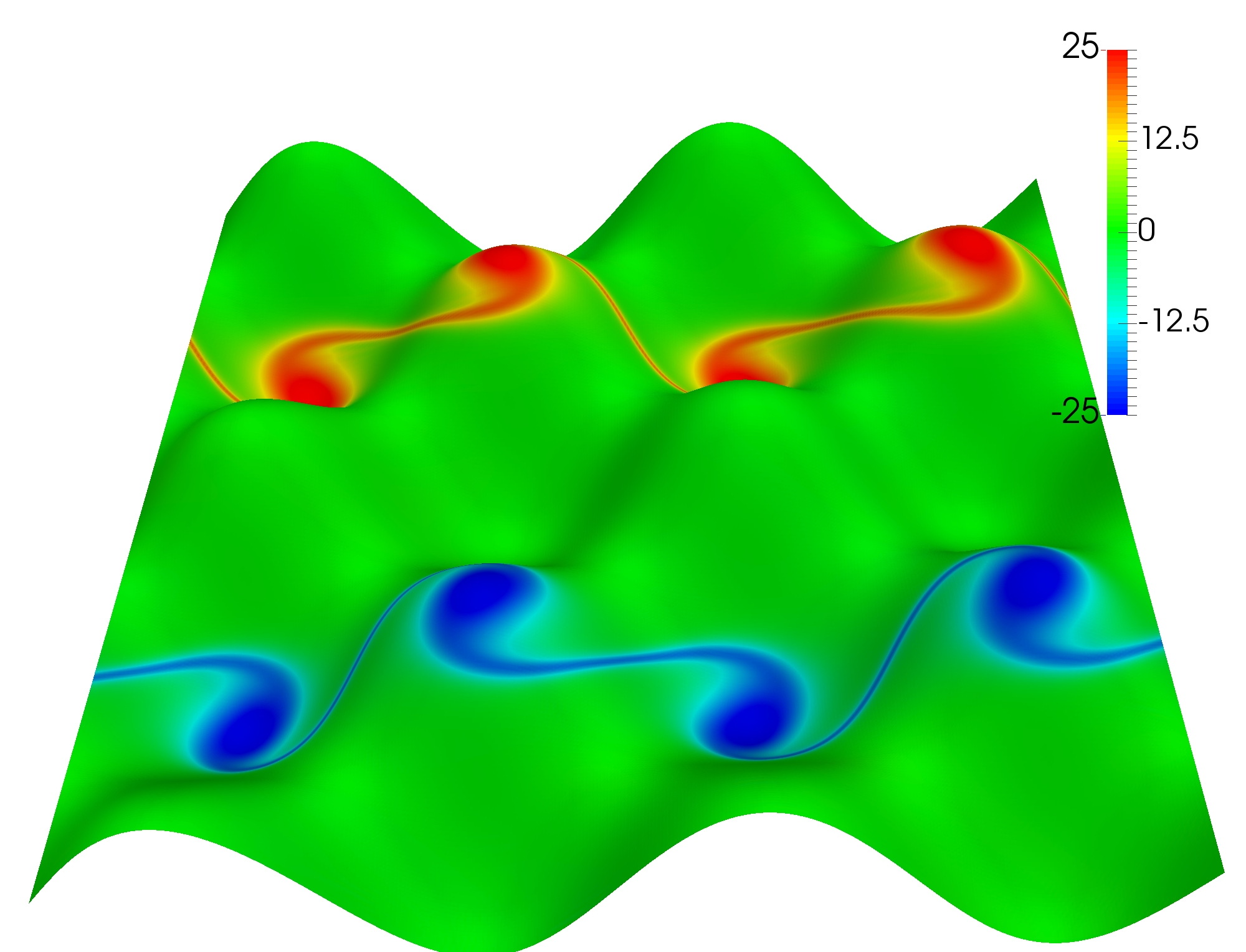} 
		\caption{}
		\label{subFig:doubleShearCurved-c}
	\end{subfigure} %
	\caption{The simulation of an inviscid incompressible flow over a curved non-Delaunay surface mesh using the circumcentric dual and the diagonal Hodge star definitions. (\subref{subFig:doubleShearCurved-a}) The curved surface showing the regions of Delaunay triangles colored in blue and the non-Delaunay triangles colored in red.  (\subref{subFig:doubleShearCurved-b}) The vorticity contour plot for the initial flow condition of a double shear layer. (\subref{subFig:doubleShearCurved-c}) The vorticity contour plot after simulation time $T = 0.28$.}
	\label{Fig:doubleShearCurved}
\end{figure}

Figure \ref{subFig:doubleShearCurved-a} shows a non-Delaunay triangulation of a sinusoidal curved surface. The surface was generated by first generating a structured-triangular mesh (consisting of isosceles right triangles) on a unit square domain and then setting the z-coordinate value for the nodes as $z= 0.1 \sin(4 \pi x) \cos (4 \pi y)$. The resulting triangulation is non-Delaunay, with almost $50 \%$ of its triangles being non-Delaunay. The Delaunay and non-Delaunay triangles are colored in blue and red in Fig. \ref{subFig:doubleShearCurved-a}. The number of elements in this curved mesh is almost equal to that for the finest Delaunay mesh. The minimum values of the primal edges length, dual edges length, primal triangles area, dual cells area are $0.003$, $2.2 \times 10^{-6}$, $4.8 \times 10^{-6}$ and $3.9 \times 10^{-6}$, respectively, which are very close to these for the finest Delaunay mesh. Such values for the curved non-Delaunay mesh, in addition to those for the sequentially subdivided mesh group, support the expectation that the distorted non-Delaunay meshes considered during this research represent relatively extreme cases for non-Delaunay triangulations. However, in practical applications, the characteristics of many of the non-Delaunay triangulations will not be very different from that for Delaunay meshes. The DEC solutions using the circumcentric dual are therefore expected to perform effectively without serious difficulty as far as the condition number of the stiffness matrices is concerned.         

The inviscid incompressible flow of a double shear layer is simulated on the curved surface with periodic boundary conditions considered on all boundaries. The initial velocity vector field is expressed as 

\begin{equation}
\label{eq:double-shear-layer}
\begin{aligned}
u_x &=\begin{cases}
\tanh ((y-0.25)/\rho), & \text{for $y \leq 0.5$},\\
\tanh ((0.75-y)/\rho), & \text{for $y > 0.5$},
\end{cases} \\
u_y &= u_z = 0.0,   
\end{aligned}
\end{equation}
with $\rho = 1/30$. Figures \ref{subFig:doubleShearCurved-b} and \ref{subFig:doubleShearCurved-c} show the vorticity contour plot at simulation times $T=0.0$ and $T=0.28$, respectively. Each shear layer developed into four vortices at the surface hills and valleys. Based on a visual metric, no difference in the solution performance can be noticed in the Delaunay versus non-Delaunay regions. Under the impression that DEC solutions using the circumcentric dual require a Delaunay mesh, one would spend considerable effort in generating a Delaunay triangulation for this curved surface. The current research shows however that a Delaunay mesh is not necessary, which significantly simplifies the mesh generation process on curved surfaces.

\section{Conclusions}
\label{sec:Conclusions}

The discrete exterior calculus solutions on non-Delaunay triangulations using circumcentric duals were numerically investigated. The present research empirically examined a common untested notion that the DEC scheme (and also the covolume method) requires a Delaunay triangulation when using the circumcentric dual. This was for long thought to be a main limitation for the DEC paradigm. The investigation was carried out through numerical convergence analysis of the 0-form Poisson equation, the 1-form Poisson equation and the incompressible Navier Stokes equations. A group of Delaunay meshes and three groups of non-Delaunay triangulations were used during the analysis. The non-Delaunay mesh groups were significantly distorted in that they had non-Delaunay triangles ratios up to $43\%$ and maximum triangles aspect ratios up to $288$. Circumcentric duals were used during all numerical experiments, where signed volumes were considered for some dual edges and some dual cell sectors in the case of non-Delaunay triangulations. Such signed volumes were used to define the diagonal Hodge star operators. 

The numerical convergence tests revealed convergent solutions for all mesh groups, with insignificant differences between the Delaunay versus the non-Delaunay groups. Furthermore, a super-convergence was observed during the incompressible Navier Stokes solutions using a group of non-Delaunay triangulations generated by sequentially subdividing a coarse Delaunay mesh. The super-convergence is attributed to coincidence of the midpoints of both the primal and dual edges that the subdivision process results in. The presented results clearly indicate that the circumcentric dual and the diagonal Hodge star definitions produce correct numerical solutions even with non-Delaunay triangulations. Such results are of practical importance since it relieves the DEC scheme from the widespread notion that it is either limited to Delaunay meshes or it requires alternative non-diagonal Hodge definitions.

The only limitation that the circumcentric duals remain to experience is the possibility of having zero volumes for some dual edges and some dual cells. The zero-volume dual edges can appear even for Delaunay triangulations, and can be avoided by minor mesh preprocessing. Such mesh preprocessing is even significantly simpler when the resulting triangulation is not restricted to be Delaunay. The zero-volume dual cells appear when the negative-volume sectors equalize the positive-volume sectors dual to a primal node. However, as pointed out earlier, significant impact for the negative-volume sectors is expected only for highly distorted triangulations, and can also be avoided through mesh processing. Again, not restricting the resulting triangulation to be Delaunay simplifies such mesh preprocessing operations significantly.     

The stiffness matrix condition number for the primal 0-form Poisson equation and the incompressible Navier Stokes equations exhibited an increase for the non-Delaunay mesh groups, in comparison to the Delaunay group. Such an increase was attributed to the decrease in both the minimum primal edges lengths and the minimum dual cells areas for the non-Delaunay groups. The significant decrease in such two geometric  attributes is not, however, a characteristic feature for non-Delaunay triangulations, and is a result of the conducted non-Delaunay mesh generation process which deliberately intended for non-Delaunay triangulations of inferior quality. In practical applications however, such geometric attributes for many non-Delaunay triangulations are not expected to be very different from corresponding Delaunay meshes. The observed behavior of both the sequentially divided non-Delaunay mesh group and the non-Delaunay curved surface triangulation support the aforementioned expectation. A more complicated mesh generation schemes would have to be considered before in order to generate a Delaunay mesh that remains Delaunay after subdivision, or to generate a Delaunay triangulation on curved surfaces. The presented research, however, relieves both the DEC paradigm and the covolume method from such longstanding unnecessary precautions.

\section*{Acknowledgments}

This research was supported by the KAUST Office of Competitive Research Funds under Award No. URF/1/1401-01-01.



\section*{References}
\bibliographystyle{elsarticle-num} 
\bibliography{DEC-References-Abbr}





\end{document}